\documentclass[10pt,a4paper]{article}

\usepackage[utf8]{inputenc}
\usepackage[T1]{fontenc} 
\usepackage[left=2cm,right=2cm,top=2cm,bottom=2cm]{geometry}

\usepackage{authblk}

\usepackage{graphicx}
\usepackage{multirow}
\usepackage{amsmath,amssymb,amsfonts}
\usepackage{amsthm}
\usepackage{mathrsfs}
\usepackage[page]{appendix}
\usepackage{textcomp}
\usepackage{manyfoot}
\usepackage[rightcaption]{sidecap}
\usepackage{physics} 
\usepackage{lineno}
\usepackage{hyperref}
\usepackage{todonotes}

\renewcommand{\epsilon}{\varepsilon}
\renewcommand{\theta}{\vartheta}
\renewcommand{\rho}{\varrho}
\renewcommand{\phi}{\varphi}

\newcommand{\bb}[1]{\mathbf{#1}}
\newcommand{\bs}[1]{\boldsymbol{#1}}

\title{Adjoint Sensitivities for the Optimization of Nonlinear Structural Dynamics via Spectral Submanifolds}

\author[1]{Matteo Pozzi}
\author[1]{Jacopo Marconi\thanks{jacopo.marconi@email.com}}
\author[2]{Shobhit Jain}
\author[3]{Mingwu Li}
\author[1]{Francesco Braghin}

\affil[1]{\footnotesize Department of Mechanical Engineering, Politecnico di Milano, Via G. La Masa 1, Milan, 20156, Italy.}
\affil[2]{Delft Institute of Applied Mathematics, TU Delft, Mekelweg 4, Delft, 2628CD, ZH, The Netherlands.}
\affil[3]{Department of Mechanics and Aerospace Engineering, Southern University of Science and Technology, Shenzhen, 518055, China.}

\date{\today} 

\begin{document}

\maketitle

\begin{abstract}
    This work presents an optimization framework for tailoring the nonlinear dynamic response of lightly damped mechanical systems using Spectral Submanifold (SSM) reduction. We derive the SSM-based backbone curve and its sensitivity with respect to parameters up to arbitrary polynomial orders, enabling efficient and accurate optimization of the nonlinear frequency-amplitude relation. We use the adjoint method to derive sensitivity expressions,  which drastically reduces the computational cost compared to direct differentiation as the number of parameters increases. An important feature of this framework is the automatic adjustment of the expansion order of SSM-based ROMs using user-defined error tolerances during the optimization process. We demonstrate the effectiveness of the approach in optimizing the nonlinear response over several numerical examples of mechanical systems. Hence, the proposed framework extends the applicability of SSM-based optimization methods to practical engineering problems, offering a robust tool for the design and optimization of nonlinear mechanical structures.
\end{abstract}

\vspace{0.5cm} 


\section{Introduction}\label{sec:introduction}

The importance of considering nonlinear effects in the design phase has been widely demonstrated by the growing number of applications where nonlinearities are not merely a side effect to be avoided but play a functional role. Notable examples include frequency division \cite{Qalandar2014frequency}, vibration mitigation \cite{BELLET20102768}, nonlinear energy sinks and targeted energy transfer \cite{Vakakis2022-ti}, micro- and nano-mechanical resonators \cite{Asadi2021strong, Li2024strain}, and micro electro-mechanical sensors \cite{ZHANG2002139, Marconi2021gyro}.

The increasing complexity of nonlinear applications has driven the advancement of analytical and computational tools in the field of Reduced Order Models (ROM, \cite{TisoMahdiabadi, touze2021model}) for analyzing the dynamic behavior of nonlinear mechanical systems. A common approach involves computing the frequency-amplitude relation of a nonlinear normal mode (NNM) through \textit{backbone curves} (see, e.g., \cite{breunung2018explicit}). These curves are typically obtained using numerical continuation techniques with the collocation method \cite{Dankowicz2013}, the Harmonic Balance Method (HBM, \cite{Krack2019}), or the shooting method \cite{Kerschen2008}. Alternatively, \textit{backbone curves} can be derived from the reduced dynamics on the Spectral Submanifold (SSM), which has proven particularly effective for high-dimensional mechanical systems \cite{Jain2022, Thurner2024}.

The variety of numerical methods, along with the ever-growing interest in nonlinear applications, has provided a great playground for developing optimization techniques aimed at tailoring the nonlinear dynamic response of mechanical systems. For instance, approaches based on the Nonlinear Normal Modes (NNMs, \cite{Kerschen2009, Haller2016-wm}) and the HBM have been used to optimize for the hardening behavior in resonators and plane frame structures \cite{DOU201663, DOU2015239, Dou2015structural}. Similarly, shape optimization has been applied to adjust eigenfrequencies and modal coupling coefficients in geometrically nonlinear MEMS gyroscopes \cite{schiwietz2024shape, schiwietz2024nonlinear}. Other strategies include nonlinear synthesis for backbone curve tailoring \cite{Detroux2021} and gradient-free optimization \cite{Li2025finite}. In the field of topology optimization \cite{Bendsoe2004}, the Equivalent Static Load method \cite{Choi2002} has been used by \cite{KIM2010660, Lee2015, Lu2021}, while \cite{Li2022eigenvalue, Dalklint2020} carried out eigenvalue optimization with frequency dependent material properties.

Despite these advancements, tailoring the nonlinear dynamic response of a system remains a significant challenge. This is mainly due to the difficulties of the numerical methods in handling high-dimensional systems with a large numbers of parameters. Additionally, these optimization approaches are highly sensitive to user-defined parameters, such as the number of harmonics in the HBM or the arc-length parameter in numerical continuation. As the system behavior evolves throughout the optimization process, these parameters may require adjustment to ensure accuracy and convergence. Moreover, inherent limitations of the numerical methods, particularly in HBM sensitivity analysis, can influence the reliability of the optimization results \cite{Saccani2022}.

Recently, Pozzi et al. \cite{Pozzi2024backbone} demonstrated the benefits of using SSM-based ROMs in optimization routines. A key advantage of this approach is the availability of analytical expressions for both the backbone curve and its sensitivity, enabling more efficient and accurate optimization. Additionally, the SSM reduction process requires only a single user-defined parameter, the expansion order, which can be automatically adjusted during the optimization based on the residual of the invariance equation. However, in this work \cite{Pozzi2024backbone}, the sensitivities were still computed via direct differentiation, which is computationally demanding as the number of optimization parameters increases. Moreover, the SSM formulation was based on the tensor notation \cite{Jain2022}, which is not the most computationally efficient choice for approximating multivariate polynomials. The second issue can be addressed by employing the multi-index notation \cite{Thurner2024}, which was proved to be more efficient than the tensor approach. This has been shown in \cite{Pozzi2025gamma}, where the authors applied topology optimization to tailor the hardening/softening dynamic response of conservative, nonlinear mechanical systems. By utilizing the adjoint sensitivity expressions derived in \cite{Li2024adjoint}, simplified to the undamped case, the approach was successfully applied to a system with 20000 design variables and around 40000 degrees of freedom. However, the explicit formulation used in \cite{Pozzi2025gamma} is limited to cubic-order computation of the Lyapunov Subcenter Manifold (LSM) in conservative systems\footnote{In the absence of damping and internal resonances, the SSM formulation reduces to that of the associated LSM.}. Furthermore, this method focuses on tailoring the backbone curve in the \textit{reduced space}, which does not allow for an optimization of  the frequency-amplitude relation in physical coordinates. Indeed, the shape of backbone-curves in physical and reduced coordinates can differ significantly, as the latter depends on the choice of parametrization \cite{Jain2022}.

To address these limitations, we focus on optimization using two-dimensional SSMs up to arbitrary polynomial orders for lightly damped, nonlinear mechanical systems. We derive expressions for adjoint sensitivities of the SSM as well as for the backbone curve in the physical space, thus allowing us to directly optimize for the physical frequency-amplitude relation. We derive these expression in the multi-index notation for a computationally efficient representation of SSMs \cite{Thurner2024}. Moreover, since the sensitivity expression is generalized for arbitrary expansion orders, we show how the accuracy of the ROM can be controlled during optimization by automatically adjusting the expansion order.

The rest of the paper is organized as follows: Sections~\ref{sec:system} and~\ref{sec:ssm_computation} introduce the mechanical system under consideration and outline the procedure for computing the SSM and the backbone curve. Section~\ref{sec:sensitivity_analysis} then presents the expressions for the backbone curve sensitivity using the adjoint method. For comparison, the sensitivity derived through direct differentiation is also provided in Appendix~\ref{app:direct_differentiation}. Finally, Section~\ref{sec:numerical_examples} showcases numerical examples to validate the proposed sensitivity formulations. Section~\ref{sec:conclusion} closes the paper with the conclusions and outlines future research directions.

\section{Mechanical system} \label{sec:system}

Consider the equation of motion of a lightly damped, autonomous, nonlinear mechanical systems:
\begin{equation} \label{eq:system_second}
    \bb{M} \ddot{\bs{x}} + \bb{C} \bs{x} + \bb{K} \bs{x} + \bs{f}(\bs{x}) = \bs{0},
\end{equation}
where $\bs{x} \in \mathbb{R}^n$ is the state vector, $n$ is the number of degrees of freedom, $\bb{M}$, $\bb{C}$, and $\bb{K}$ are the mass, damping, and stiffness matrices (all belonging to $\mathbb{R}^{n \times n}$), and $\bs{f} \in \mathbb{R}^n$ is the vector of displacement-dependent nonlinear force. In this work, $\bs{f}$ contains the contribution of quadratic and cubic nonlinearities in the displacement
\begin{equation} \label{eq:nonlinear_force}
    \bs{f}(\bs{x}) = \bs{f}_2(\bs{x}, \bs{x}) + \bs{f}_3(\bs{x}, \bs{x}, \bs{x}),
\end{equation}
and represents the contribution of the internal nonlinear elastic forces. These are efficiently managed using the tensor notation \cite{Marconi2020} as
\begin{equation} \label{eq:nonlinear_force_tensor}
    f^i = T_2^{ijk} x^j x^k + T_3^{ijkl} x^j x^k x^l,
\end{equation}
where $f^i$ is the $i^{th}$ element of $\bs{f}$, and the tensors $T_2 \in \mathbb{R}^{n \times n \times n}$ and $T_3 \in \mathbb{R}^{n \times n \times n \times n}$ contain the coefficients of the quadratic and cubic nonlinearities. Notice that we denote the components of a vector or a tensor using superscripts (not to be confused with exponents).

The damping matrix $\bb{C}$ is computed using the Rayleigh \cite{oda2008} damping formulation:
\begin{equation} \label{eq:damping_matrix}
    \bb{C} = \alpha_R \bb{M} + \beta_R \bb{K},
\end{equation}
where $\alpha_R$ and $\beta_R$ are non-negative coefficients (see Appendix~\ref{app:damping} for more details on the choice of the damping coefficients).

The second-order system~\eqref{eq:system_second} is transformed into the first-order form:
\begin{equation} \label{eq:system_first}
    \bb{B} \dot{\bs{z}} = \bb{A} \bs{z} + \bs{F}(\bs{z}),
\end{equation}
where
\begin{equation}
\begin{aligned}
    \bs{z} = 
    \begin{bmatrix}
        \bs{x} \\
        \dot{\bs{x}}
    \end{bmatrix},\quad 
    \bb{B} = 
    \begin{bmatrix}
        \bb{C} & \bb{M} \\
        \bb{M} & \bb{0}
    \end{bmatrix},\quad
    \bb{A} =
    \begin{bmatrix}
        -\bb{K} & \bb{0} \\
        \bb{0} & \bb{M}
    \end{bmatrix}
    ,\quad \bs{F}(\bs{z}) = 
    \begin{bmatrix}
        -\bs{f}(\bs{x}) \\
        \bs{0}
    \end{bmatrix}.
\end{aligned}
\end{equation}

Assume that the autonomous part of system~\eqref{eq:system_first} has a stable and hyperbolic fixed point at $\bs{z} = \bs{0}$. The eigenvalues $\lambda_i \in \mathbb{C}$ of such system are computed by solving the generalized right and left eigenvalue problems:
\begin{equation} \label{eq:gep_first}
    \left( \bb{A} - \lambda_i \bb{B} \right) \bs{v}_i = \bs{0} ,\quad
    \bs{u}_i^* \left( \bb{A} - \lambda_i \bb{B} \right) = \bs{0},
\end{equation}
where $(\bullet)^*$ represents the conjugate transpose operator, and $\bs{v}_i \in \mathbb{C}^{2n}$ and $\bs{u}_i \in \mathbb{C}^{2n}$ are the right and left eigenvectors associated to $\lambda_i$.

The eigenvalues and eigenvectors of the first order system (computed from Eq.\eqref{eq:gep_first}) can be related to those of the respective undamped second order system as:
\begin{equation} \label{eq:eigenvalues}
    \lambda_1 = \bar{\lambda}_2 = \lambda = -\xi\omega + \mathrm{i}\omega\sqrt{1 - \xi^2} = \alpha + \mathrm{i} \omega_d
\end{equation}
and
\begin{equation} \label{eq:eigenvectors}
    \bs{v}_1 = \bar{\bs{v}}_2 =
    \begin{bmatrix}
        \bs{\phi} \\ \lambda \bs{\phi}
    \end{bmatrix},\quad
    \bs{u}_1 = \bar{\bs{u}}_2 =
    \begin{bmatrix}
        \bs{\phi} \\ \bar{\lambda} \bs{\phi}
    \end{bmatrix},
\end{equation}
where $\bs{\phi} \in \mathbb{R}^n$ and $\omega \in \mathbb{R}$ are obtained from the generalized eigenvalue problem
\begin{equation} \label{eq:gep}
    \left(\bb{K} - \omega^2 \bb{M}\right) \bs{\phi} = \bs{0},
\end{equation}
where the mode shapes $\bs{\phi}$ are assumed to be mass normalized, yielding
\begin{equation} \label{eq:mass_normalization}
    \bs{\phi}^T \bb{M} \bs{\phi} = 1 , \quad \bs{\phi}^T \bb{K} \bs{\phi} = \omega^2. \\
\end{equation}

Using the definition of the damping matrix and the mass normalization condition (Eqs.~\eqref{eq:damping_matrix} and~\eqref{eq:mass_normalization}), the damping ratio of the system is defined as
\begin{equation} \label{eq:damping_ratio}
    \xi = \frac{\bs{\phi}^T \bb{C} \bs{\phi}}{2\omega} = \frac{\alpha_R + \beta_R \omega^2}{2\omega}.
\end{equation}

For lightly damped systems, we require that $\xi \ll 1$.

\section{SSM computation} \label{sec:ssm_computation}

We provide here the basic steps used to compute the \textit{two-dimensional} SSM for an autonomous, damped, nonlinear mechanical system (see \cite{Jain2022, Thurner2024} for a complete analysis).

Call $E$ the two-dimensional master subspace spanned by $\bs{v}_1$ and $\bs{v}_2$. Let $\bs{R}:\mathbb{C}^2 \to \mathbb{C}$ be a two-dimensional parametrization of the reduced dynamics, such that
\begin{equation} \label{eq:reduced_dynamics}
    \dot{\bs{p}} = \bs{R}(\bs{p}),
\end{equation}
where $\bs{p} \in \mathbb{C}^2$ is the vector of reduced coordinates. Using $\bs{p}$, let $\bs{W}:\mathbb{C}^2 \to \mathbb{R}^{2n}$ be a parametrization of the SSM attached to the master subspace $E$. Then, the state vector $\bs{z}$ can be expressed as
\begin{equation} \label{eq:manifold}
    \bs{z} = \bs{W}(\bs{p}) = 
    \begin{bmatrix}
        \bs{w}(\bs{p}) \\ \dot{\bs{w}}(\bs{p})
    \end{bmatrix},
\end{equation}
where $\bs{w}(\bs{p})$ and $\dot{\bs{w}}(\bs{p})$ are the parametrizations for the displacement $\bs{x}$ and for the velocity $\dot{\bs{x}}$, respectively.

Substituting the parametrizations~\eqref{eq:reduced_dynamics} and~\eqref{eq:manifold} into Eq.~\eqref{eq:system_first}, we obtain the \textit{autonomous invariance equation}:
\begin{equation} \label{eq:invariance}
    \bb{B} \pdv{\bs{W}(\bs{p})}{\bs{p}} \bs{R}(\bs{p}) = \bb{A} \bs{W}(\bs{p}) + \bs{F}(\bs{W}(\bs{p})).
\end{equation}

To numerically solve Eq.~\eqref{eq:invariance}, the parametrizations~\eqref{eq:reduced_dynamics} and~\eqref{eq:manifold} are Taylor expanded up to a given order:
\begin{align}
    \bs{R}(\bs{p}) \approx \sum_{\bs{m} \in \mathbb{N}^2} \bs{R}_{\bs{m}} \bs{p}^{\bs{m}}, \quad\quad
    \bs{w}(\bs{p}) \approx \sum_{\bs{m} \in \mathbb{N}^2} \bs{w}_{\bs{m}} \bs{p}^{\bs{m}}, \quad\quad
    \dot{\bs{w}}(\bs{p}) \approx \sum_{\bs{m} \in \mathbb{N}^2} \dot{\bs{w}}_{\bs{m}} \bs{p}^{\bs{m}},
\end{align}
where the multi-index notation is used to express the multivariate polynomials (see Appendix~\ref{app:multi_index}). The terms $\bs{R}_{\bs{m}} = \{R_{\bs{m}}^1, R_{\bs{m}}^2\}^T \in \mathbb{C}^2$ are the coefficients of the reduced dynamics, while $\bs{w}_{\bs{m}} \in \mathbb{C}^n$ and $\dot{\bs{w}}_{\bs{m}} \in \mathbb{C}^n$ are the manifold coefficients.

Following the approach described in \cite{Thurner2024}, it is possible to split the \textit{autonomous invariance equation} into many \textit{cohomological equations}, one for each multi-index $\bs{m}$:
\begin{align}
    \bb{L}_{\bs{m}} \bs{w}_{\bs{m}} &= \bs{h}_{\bs{m}} \label{eq:cohomological_equation} \\
    \bb{L}_{\bs{m}} &= \bb{K} + \Lambda_{\bs{m}} \bb{C} + \Lambda_{\bs{m}}^2 \bb{M} \\
    \Lambda_m &= \bs{m} \cdot \bs{\Lambda} \\
    \bs{h}_{\bs{m}} &= \bs{D}_{\bs{m}} \bs{R}_{\bs{m}} + \bs{C}_{\bs{m}},
\end{align}
where $\bs{\Lambda} = \{ \lambda_1, \lambda_2 \}^T$, while $\bs{C}_{\bs{m}} \in \mathbb{C}^n$ and $\bs{D}_{\bs{m}} \in \mathbb{C}^{n \times 2}$ are defined in Appendix~\ref{app:ssm_steps}. $\bs{C}_{\bs{m}}$, in particular, depends on the parametrization coefficients at orders lower than $m$.

Equation~\eqref{eq:cohomological_equation} is solved up to any arbitrary expansion order for the manifold coefficients $\bs{w}_{\bs{m}}$, starting from the first. At each order $m$, there are $m+1$ distinct multi-indices. Details on the computation of $\bs{R}_{\bs{m}}$ and $\dot{\bs{w}}_{\bs{m}}$ are given in Appendix~\ref{app:ssm_steps}. Here, we just mention that the coefficients $\bs{R}_{\bs{m}}$ are obtained following the normal-form style parametrization \cite{Thurner2024}, for which the coefficients are nonzero only if the near-resonance condition is satisfied (see Appendix~\ref{app:near_resonance}). This allows to write an \textit{analytical expression} for the backbone curve \cite{Pozzi2024backbone}.

\subsection{Backbone curve} \label{sec:backbone}

Using the polar version of the reduced coordinates $\bs{p}(\theta) = \rho \{ e^{\mathrm{i}\theta}, e^{-\mathrm{i}\theta} \}^T$, it is possible to write an analytical expression for the backbone curve in the reduced space \cite{Pozzi2024backbone}:
\begin{equation} \label{eq:Omega_compact}
    \Omega(\rho) = \Im{\lambda} + \sum_{m > 1} \Im{R_{\bs{m}}^1} \rho^{m - 1},
\end{equation}
where the response frequency $\Omega$ is function of the reduced amplitude $\rho$. According to the normal-form style parametrization \cite{Thurner2024}, the coefficients $R_{\bs{m}}^j$ are non-zero if the multi-index $\bs{m}$ satisfies the near-resonance condition (Appendix~\ref{app:near_resonance}). From this condition, it follows that the coefficients associated to even orders are always zero, meaning that the summation in Eq.~\eqref{eq:Omega_compact} can be performed only on odd orders.

For the purposes of computing the sensitivities, the expression for $\Omega$ can be rewritten as
\begin{equation} \label{eq:Omega}
        \Omega(\rho) = \frac{1}{2} \mathrm{i} \left( \lambda_2 - \lambda_1 \right) + \frac{1}{2} \mathrm{i} \sum_{m > 1} \left( R_{\bs{m}}^2 - R_{\bs{m}}^1 \right) \rho^{m - 1}.
\end{equation}

Let us write a measure for the physical amplitude. Call $x^i$ the $i^{th}$ element of the vector $\bs{x}$, corresponding to the displacement of the $i^{th}$ degree of freedom. Similarly, let us call $w_{\bs{m}}^i$ the $i^{th}$ component of the vector $\bs{w}_{\bs{m}}$. This value can be computed using the manifold parametrization (Eq.~\eqref{eq:manifold}):
\begin{equation}
    x^i(\rho, \theta) = \sum_{\bs{m}} w^i_{\bs{m}} \bs{p}^{\bs{m}}(\varrho,\theta),
\end{equation}
where $x^i$ is a function of the reduced amplitude $\rho$ and of the phase $\theta$. To remove the dependence on $\theta$, we define a set of $\theta$ values as
\begin{equation} \label{eq:theta_values}
    \theta_k = k \frac{2\pi}{N_{\theta}} \quad \forall k = 1, \dots, N_{\theta}.
\end{equation}

Then, we compute the Root Mean Square (RMS) of $x^i$ over the values of $\theta_k$, leading to
\begin{equation} \label{eq:x_rms}
    x(\rho) = \sqrt{\frac{1}{N_{\theta}} \sum_{k = 1}^{N_{\theta}}\left( x^i(\rho, \theta_k)\right)^2}.
\end{equation}

In an optimization framework, this relation can be used to find the reduced amplitude $\rho$ corresponding to a target physical amplitude $x=x_0$. The obtained $\rho$ can be then used to evaluate the corresponding frequency $\Omega$ on the backbone from Eq.~\eqref{eq:Omega}. This way, the frequency-amplitude relation can be tracked in the \textit{physical} space.

\section{Sensitivity analysis} \label{sec:sensitivity_analysis}

As done in \cite{Pozzi2024backbone}, we propose an optimization procedure to tailor the backbone curve by defining target points in terms of response frequency $\Omega$ and physical amplitude $x$. An example of optimization problem can be stated as
\begin{equation} \label{opt:optimization}
    \begin{aligned}
        \min_{\bs{\mu}} \quad& J \\
        \mathrm{s.t.} \quad& \Omega(\rho(x_{j})) - \Omega_j = 0, \quad \forall j \\
        & \bs{\mu}_L < \bs{\mu} < \bs{\mu}_U,
    \end{aligned}
\end{equation}
where $J$ is the objective function, $\Omega_j$ is the $j^{th}$ target frequency corresponding to the target amplitude $x_j$, and $\bs{\mu}$ is the vector of design variables with lower and upper bounds $\bs{\mu}_L$ and $\bs{\mu}_U$.

To solve the optimization problem~\eqref{opt:optimization} using a gradient-based approach, the sensitivities of both the objective function $J$ and the constraints must be computed. In particular, the sensitivity of the response frequency $\Omega$ depends on the SSM coefficients. A straightforward approach to obtain this sensitivity is the direct differentiation method, which involves applying the chain rule to differentiate all equations contributing to the computation of $\Omega$. The detailed steps of this procedure are provided in Appendix~\ref{app:direct_differentiation}. However, the computational cost of direct differentiation scales with the number of design variables, which can limit the size of the optimization problem.

To efficiently handle optimization problems with a large number of design variables, we instead compute the sensitivity of $\Omega$ using the adjoint method \cite{Manzoni2021optimal}, in which the complex derivatives are managed through \textit{Wirtinger calculus} \cite{Wirtinger1927, Manzoni2021optimal}.

\subsection{The adjoint method} \label{sec:adjoint_sensitivity}

The first step of the \textit{adjoint method} is the identification of the \textit{state variables} and the corresponding \textit{state functions}. In this case, the \textit{state variables} are the reduced amplitude $\rho$, the manifold coefficients $\bs{w}_{\bs{m}}$ for $m > 1$, the mode shape $\bs{\phi}$, and the natural frequency $\omega$. The corresponding \textit{state functions} are the RMS of the physical amplitude (Eq.~\eqref{eq:x_rms}), the cohomological equations (Eq.~\eqref{eq:cohomological}), the generalized eigenvalue problem (Eq.~\eqref{eq:gep}), and the mass normalization condition (Eq.~\eqref{eq:mass_normalization}).

Next, the \textit{Lagrangian function} is assembled as follows:
\begin{equation} \label{eq:lagrangian}
\begin{aligned}
    \mathcal{L} &= \Omega + \lambda_{\rho} \left( x - x_0 \right) + \sum_{m > 1} \bs{\lambda}_{\bs{m}}^T \left( \bb{L}_{\bs{m}} \bs{w}_{\bs{m}} - \bs{h}_{\bs{m}} \right) + \bs{\lambda}_{\bs{\phi}}^T \left( \bb{K} - \omega^2 \bb{M} \right) \bs{\phi} + \lambda_{\omega} \left( \bs{\phi}^T \bb{M} \bs{\phi} - 1 \right),
\end{aligned}
\end{equation}
where $\lambda_{\rho}$, $\bs{\lambda}_{\bs{m}}$, $\bs{\lambda}_{\bs{\phi}}$, and $\lambda_{\omega}$ are the \textit{adjoint variables}. According to \textit{Wirtinger calculus} \cite{Wirtinger1927, Manzoni2021optimal}, only the real part of the cohomological equations should be used, as the Lagrangian must be a real-valued function. However, in this case, it can be shown that the sum of all cohomological equations inherently results in a real-valued function. This follows from the \textit{symmetric multi-index} property (Appendix~\ref{app:multi_index}), which ensures that the cohomological equations appear in complex-conjugate pairs.

The \textit{adjoint variables} are computed solving the \textit{adjoint equations}, which are obtained by taking the partial derivatives of $\mathcal{L}$ with respect to the \textit{state variables}:
\begin{align}
    \pdv{\mathcal{L}}{\rho} &= \pdv{\Omega}{\rho} + \lambda_{\rho} \pdv{x}{\rho} = 0 \label{eq:adjoint_equation_rho} \\
    \begin{split} \label{eq:adjoint_equation_w}
        \pdv{\mathcal{L}}{\bs{w}_{\bs{m}}} &= \pdv{\Omega}{\bs{w}_{\bs{m}}} + \lambda_{\rho} \pdv{x}{\bs{w}_{\bs{m}}} + \bs{\lambda}_{\bs{m}}^T \bb{L}_{\bs{m}} - \sum_{q > m} \bs{\lambda}_{\bs{q}}^T \pdv{\bs{h}_{\bs{q}}}{\bs{w}_{\bs{m}}} = \bs{0}^T
    \end{split} \\
    \begin{split} \label{eq:adjoint_equation_phi}
        \pdv{\mathcal{L}}{\bs{\phi}} &= \pdv{\Omega}{\bs{\phi}} + \lambda_{\rho} \pdv{x}{\bs{\phi}} - \sum_{m > 1} \bs{\lambda}_{\bs{m}}^T \pdv{\bs{h}_{\bs{m}}}{\bs{\phi}} + \bs{\lambda}_{\bs{\phi}}^T \left( \bb{K} - \omega^2 \bb{M} \right) + 2\lambda_{\omega} \bs{\phi}^T \bb{M} = \bs{0}^T
    \end{split} \\
    \begin{split} \label{eq:adjoint_equation_omega}
        \pdv{\mathcal{L}}{\omega} &= \pdv{\Omega}{\omega} + \sum_{m > 1} \bs{\lambda}_{\bs{m}}^T \left( \pdv{\bb{L}_{\bs{m}}}{\omega} \bs{w}_{\bs{m}} - \pdv{\bs{h}_{\bs{m}}}{\omega} \right) - 2\omega \bs{\lambda}_{\bs{\phi}}^T \bb{M} \bs{\phi} = 0,
    \end{split}
\end{align}
where the equations for $\bs{\lambda}_{\bs{m}}$ are computed in reverse order starting from the highest one, and the last two equations are solved together. All the partial derivatives involved in Eqs.~\eqref{eq:adjoint_equation_rho},~\eqref{eq:adjoint_equation_w},~\eqref{eq:adjoint_equation_phi}, and~\eqref{eq:adjoint_equation_omega} are given in Appendix~\ref{app:adjoint}.

Using the adjoint variables, the sensitivity of $\Omega$ is obtained by taking the partial derivative of $\mathcal{L}$ with respect to all the design variables $\mu$:
\begin{equation} \label{eq:lagrangian_mu}
\begin{aligned}
    \pdv{\mathcal{L}}{\mu} &= \pdv{\Omega}{\mu} + \sum_{m > 1} \bs{\lambda}_{\bs{m}}^T \left( \pdv{\bb{L}_{\bs{m}}}{\mu} \bs{w}_{\bs{m}} - \pdv{\bs{h}_{\bs{m}}}{\mu} \right) + \bs{\lambda}_{\bs{\phi}}^T \left( \pdv{\bb{K}}{\mu} - \omega^2 \pdv{\bb{M}}{\mu} \right) \bs{\phi} + \lambda_{\omega} \bs{\phi}^T \pdv{\bb{M}}{\mu} \bs{\phi},
\end{aligned}
\end{equation}
where the partial derivatives with respect to $\mu$ depend on how the system is parametrized, and they are usually straightforward to compute.

It is important to point out that the \textit{symmetric multi-index} property (Appendix~\ref{app:multi_index}) also holds for the adjoint variables $\bs{\lambda}_{\bs{m}}$. In this way, we can avoid solving for nearly half of the adjoint equations associated to $\bs{\lambda}_{\bs{m}}$.

\begin{figure}[t]
    \centering
    \includegraphics[width=1\linewidth]{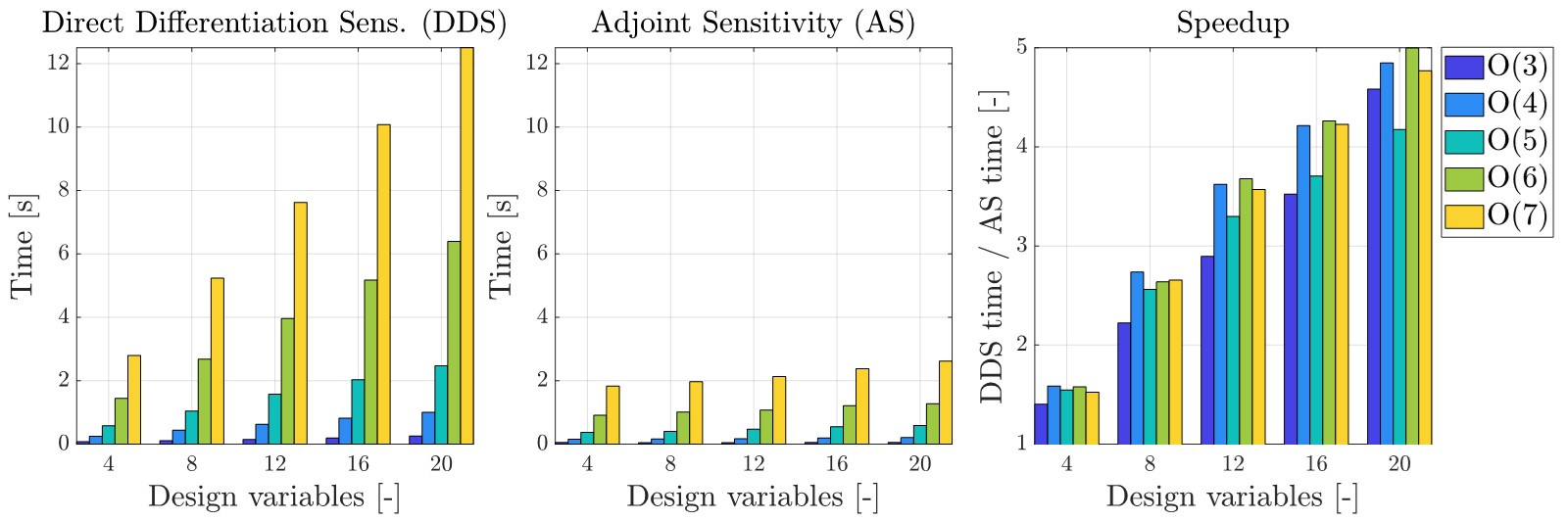}
    \caption{Computational times involved in the sensitivity analysis for a system with 501 degrees of freedom. The direct differentiation (left) and the adjoint method (center) are compared for different SSM expansion orders and number of design variables. On the right, the speedup is shown.}
    \label{fig:time_sens}
\end{figure}

Figure~\ref{fig:time_sens} illustrates the computational time involved in the sensitivity analysis. Specifically, the adjoint formulation is compared to the direct differentiation approach (outlined in Appendix~\ref{app:direct_differentiation}) for a system with 501 degrees of freedom. The comparison is made for different SSM expansion orders and varying numbers of design variables. All computations were performed on a Windows laptop with an Intel Core i7-1255U CPU @ 1.70 GHz and 16.0 GB of RAM.

As expected, direct differentiation scales linearly with the number of design variables, while the adjoint method remains independent\footnote{The adjoint sensitivity (Eq.~\eqref{eq:lagrangian_mu}) must be evaluated for each design variable, introducing a computational cost that scales linearly with their number. However, this overhead remains negligible compared to the linear scaling of direct differentiation.} of them, making it a more efficient approach for computing sensitivities.

\section{Numerical examples} \label{sec:numerical_examples}

In this Section, three numerical examples are provided to highlight the advantages of the proposed adjoint sensitivity formulation. With the first one, we validate the adjoint sensitivity by approximating the perturbed backbone curve of two coupled nonlinear oscillators. Then, we use the sensitivity expression in a gradient-based optimization loop for tailoring the backbone curve of nonlinear mechanical systems.

In all the optimization examples, the Modal Assurance Criterion (MAC) is used to identify and track the target mode shape throughout the optimization. This is required as the order of the modes may change during the optimization, thus leading to convergence issues if not properly addressed. See, for instance, \cite{Pozzi2023temperature, Pozzi2024backbone} for more details about the MAC.

Another important detail is the accuracy of the SSM during the optimization. The invariance equation is used to define an error measure $\epsilon$ that quantifies the accuracy of the SSM. By setting an error tolerance $\epsilon_{tol}$ it is possible to automatically adjust the expansion order to keep the error measure below the defined tolerance. As done in \cite{Pozzi2024backbone}, the error measure is computed according to the residual of the invariance equation \eqref{eq:invariance}.

All the examples are computed in YetAnotherFEcode v1.4.0 \cite{yafec} using Matlab 2023b.

\subsection{Perturbed backbone of two coupled nonlinear oscillators}

Consider the equations of motion of two coupled oscillators, with the first mass attached to the ground and to the second mass through nonlinear springs:
\begin{equation}
    \begin{split}
        m\ddot x_1 + c \dot x_1 + k(2x_1 - x_2) + k_2 x_1^2 + k_3x_1^3 + k_2 (x_1 - x_2)^2 + k_3 (x_1 - x_2)^3 = 0, \\
        m\ddot x_2 + c \dot x_2 + k(x_2 - x_1) + k_2 (x_2 - x_1)^2 + k_3 (x_2 - x_1)^3 = 0,
    \end{split}
\end{equation}
where $m = 1$, $k = 1$, $k_2 = 0.5$, and $k_3 = 0.2$. The damping is computed according to the Rayleigh damping with $\alpha_R = 0$ and $\beta_R = 0.1$, leading to a damping coefficient $c = \beta_R k = 0.1$. We focus on the displacement of the mass at the free end ($x_2$). The target mode shape is the first one, associated with eigenvalues $\lambda = -0.0191 \pm \mathrm{i} 0.6177$. The SSM is expanded up to the $5^{th}$ order.

\begin{figure*}[t]
    \centering
    \includegraphics[width=0.9\linewidth]{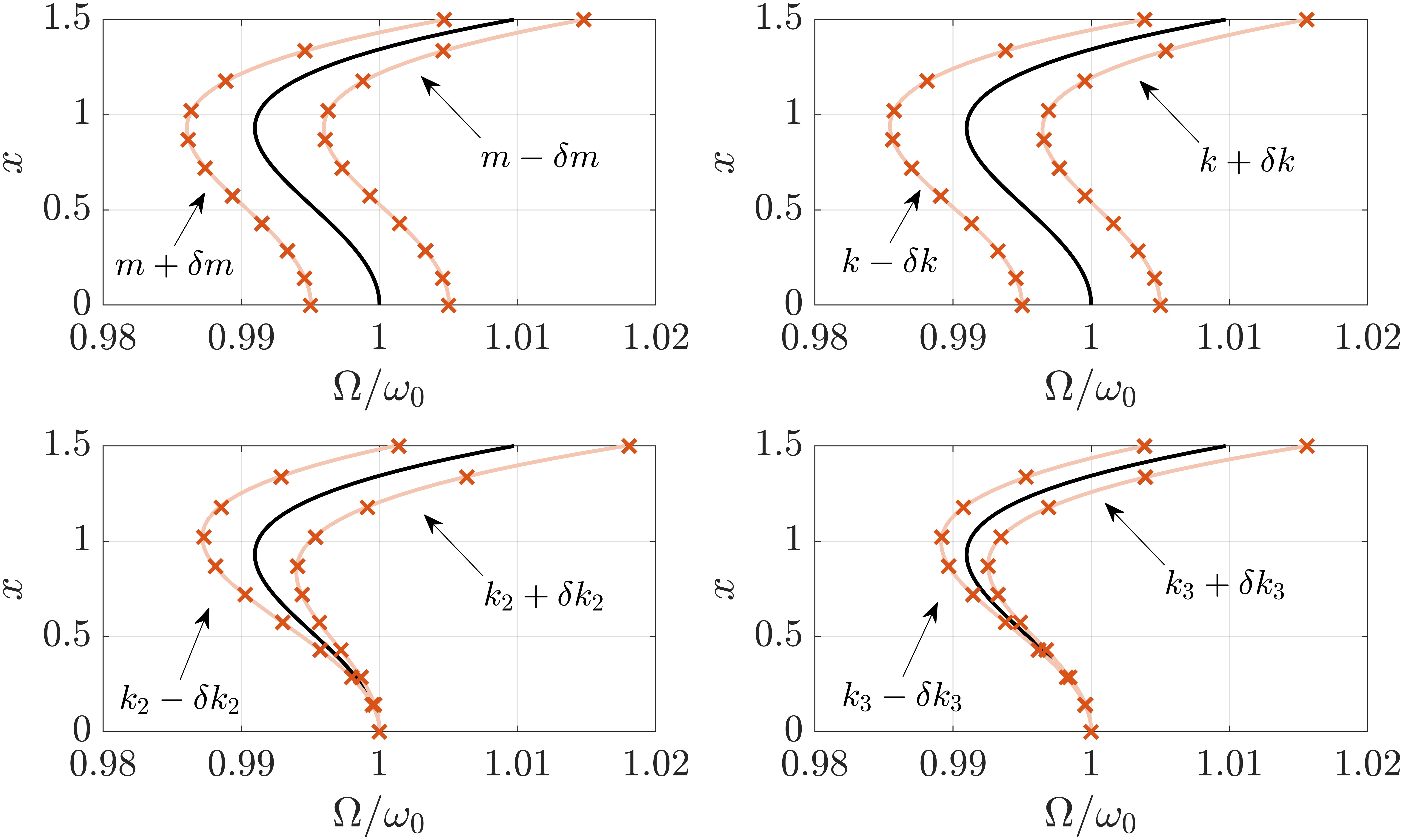}
    \caption{Nominal (black lines) and perturbed (colored lines) backbone curves for the coupled oscillators. $x$ is the RMS amplitude of $x_2$. The cross markers represent the predicted points obtained using the sensitivity values as in Eq.~\eqref{eq:backbone_prediction}. In this example, the parameter perturbations are $\delta m = 0.01\,m$, $\delta k = 0.01\,k$, $\delta k_2 = 0.03\,k_2$, and $\delta k_3 = 0.03\,k_3$.}
    \label{fig:chain_backbone}
\end{figure*}

First, the backbone curve of the system in nominal condition is computed (black lines in Fig.~\ref{fig:chain_backbone}). Then, the sensitivities of the backbone curve with respect to the system parameters $\mu$ are computed using the adjoint method. These are used to make a prediction (cross markers in Fig.~\ref{fig:chain_backbone}) of the backbone curve of the perturbed system as
\begin{equation} \label{eq:backbone_prediction}
    \Omega(\mu + \delta\mu) = \Omega(\mu) + \pdv{\Omega}{\mu} \delta\mu.
\end{equation}

The predictions are then validated by actually computing the backbone curves of the perturbed system (colored lines in Fig.~\ref{fig:chain_backbone}). The natural frequency only changes due to perturbations in $m$ and $k$, while perturbations in $k_2$ and $k_3$ become relevant at higher amplitudes, thus affecting the hardening/softening dynamic behavior of the response.

Moreover, the results also show that the sensitivity-based predictions can be reliably used to study a system around its nominal parameters, making uncertainty quantification analysis extremely efficient \cite{morsy2025}.

\subsection{Optimization of a von K\'arm\'an beam} \label{sec:von_karman}

Consider a geometrically nonlinear, undamped, clamp-clamp von K\'arm\'an beam. This example provides a basis for comparing the proposed optimization approach, which uses multi-index notation to parametrize the manifold, with the method presented in \cite{Pozzi2024backbone}, which employs tensorial notation. In both cases, the direct differentiation sensitivities, provided in \cite{Pozzi2024backbone} for the tensorial notation and in Appendix~\ref{app:direct_differentiation} for the multi-index
notation, are used.

The structure is described by a finite element model consisting of 10 von K\'arm\'an elements (27 degrees of freedom). As design variables, we take the beam thickness $h$, the length $L$, and the amplitudes $A_1$ and $A_2$ which define the shape of the beam as
\begin{equation*}
    y = A_1\sin(\pi x /L)+A_2\sin(2\pi x/L),
\end{equation*}
being $x\in[0,L]$ and $y$ the nodal coordinates of the beam. The material is characterized by Young's modulus equal to $90\,\mathrm{GPa}$, Poisson's ratio of $0.3$, and density corresponding to $7850\,\mathrm{kg/m^3}$.

\begin{figure}[t]
    \centering
    \includegraphics[width=0.8\linewidth]{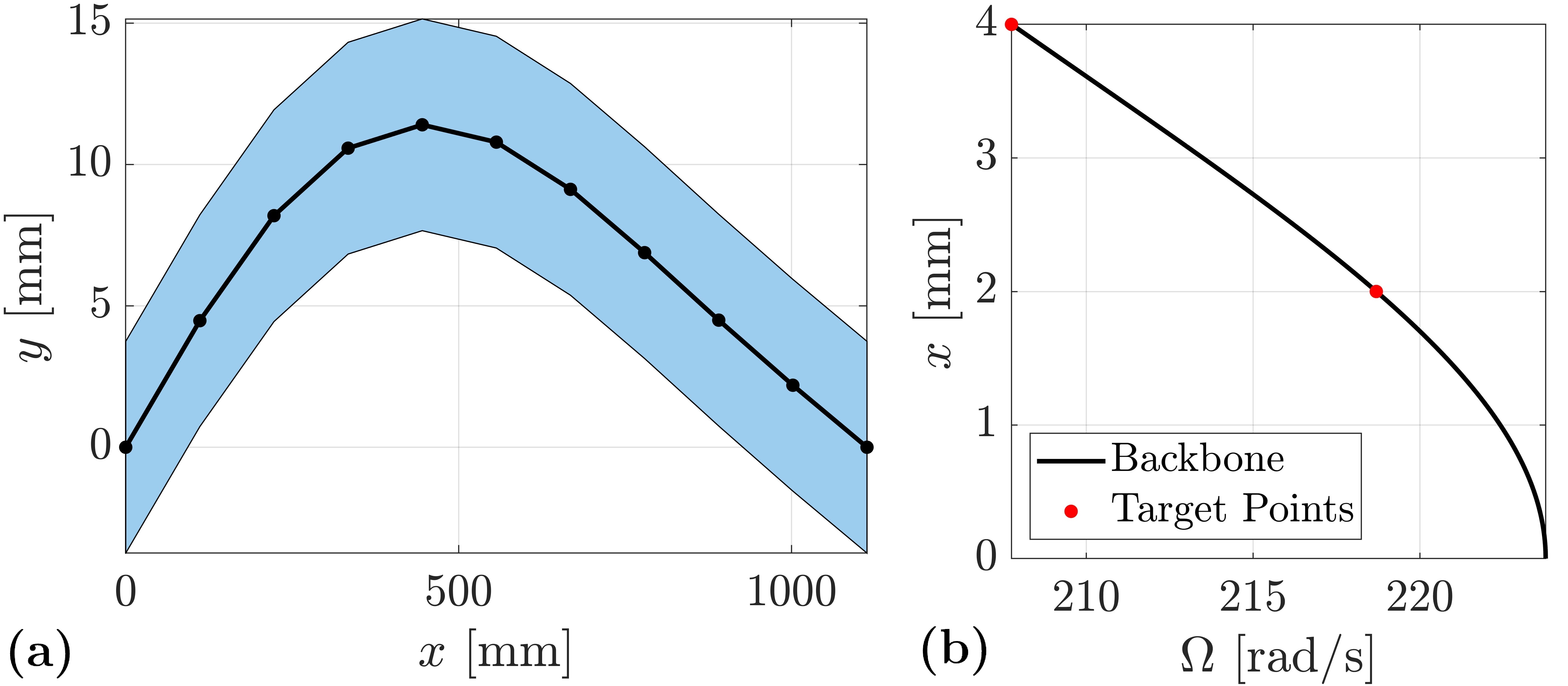}
    \caption{ Optimal solution obtained solving problem~\eqref{opt:von_karman}. (a) Optimized beam shape. (b) Backbone curve corresponding to the optimal solution. The same results were obtained in \cite{Pozzi2024backbone} in 4.8 min, while here the computations took only 15 s.}
    \label{fig:von_karman_solution}
\end{figure}

Initial values and bounds for the design variables are reported in Table~\ref{tab:von_karman_initial}. The optimization problem is stated as follows:
\begin{equation} \label{opt:von_karman}
    \begin{aligned}
        \min_{\bs{\mu}} \quad& A_2 L \\
        \mathrm{s.t.} \quad& \Omega \left( \rho \left(v=0.2\,h_0\right)\right) - \omega_0 = 0 \\
        & \Omega \left(\rho\left(v=0.4\,h_0\right)\right) - 0.95\,\omega_0 = 0 \\
        & \bs{\mu}_L \leq \bs{\mu} \leq \bs{\mu}_U,
    \end{aligned}
\end{equation}
where $h_0$ is the initial thickness, and $\omega_0$ is the first eigenfrequency of the system (evaluated using the initial parameters). $v$ is the y-displacement of the center of the beam, which is also the DOF we use to compute the backbone. In Fig.~\ref{fig:von_karman_solution}, $x$ is the RMS amplitude of $v$.

\begin{table}[t]
    \centering
    \begin{tabular}{l l r r r r}
        \hline
        \multicolumn{2}{l}{Design variable} & Initial value & Lower bound & Upper bound & Optimal value \\
        \hline
        $A_1$ & $[\mathrm{mm}]$ & 0    & 0    & 20   & 10.7   \\
        $A_2$ & $[\mathrm{mm}]$ & 0    & 0    & 20   & 1.94   \\
        $h$   & $[\mathrm{mm}]$ & 10   & 1    & 100  & 7.49   \\
        $L$   & $[\mathrm{mm}]$ & 1000 & 500  & 1500 & 1113.6 \\
        \hline
    \end{tabular}
    \caption{Initial values, lower bounds, and upper bounds of the design variables for the von K\'arm\'an beam example.}
    \label{tab:von_karman_initial}
\end{table}

Using the approach presented in \cite{Pozzi2024backbone}, the optimization reaches convergence in 11 iterations. The overall computational time was about 4.8 minutes on a Windows laptop equipped with Intel Core i7-1255U CPU @ 1.70 GHz and 16.0 GB RAM. On the same machine, solving the same problem with the direct differentiation sensitivity (Section~\ref{sec:sensitivity_analysis}) yields identical results and optimization history (Fig~\ref{fig:von_karman_solution}), but the process only takes around 15 seconds. This is because the multi-index notation is much more efficient than the tensorial one.

\subsection{Optimization of a MEMS gyroscope} \label{sec:mems_gyro}

\begin{figure}[t]
    \centering
    \includegraphics[width=.9\linewidth]{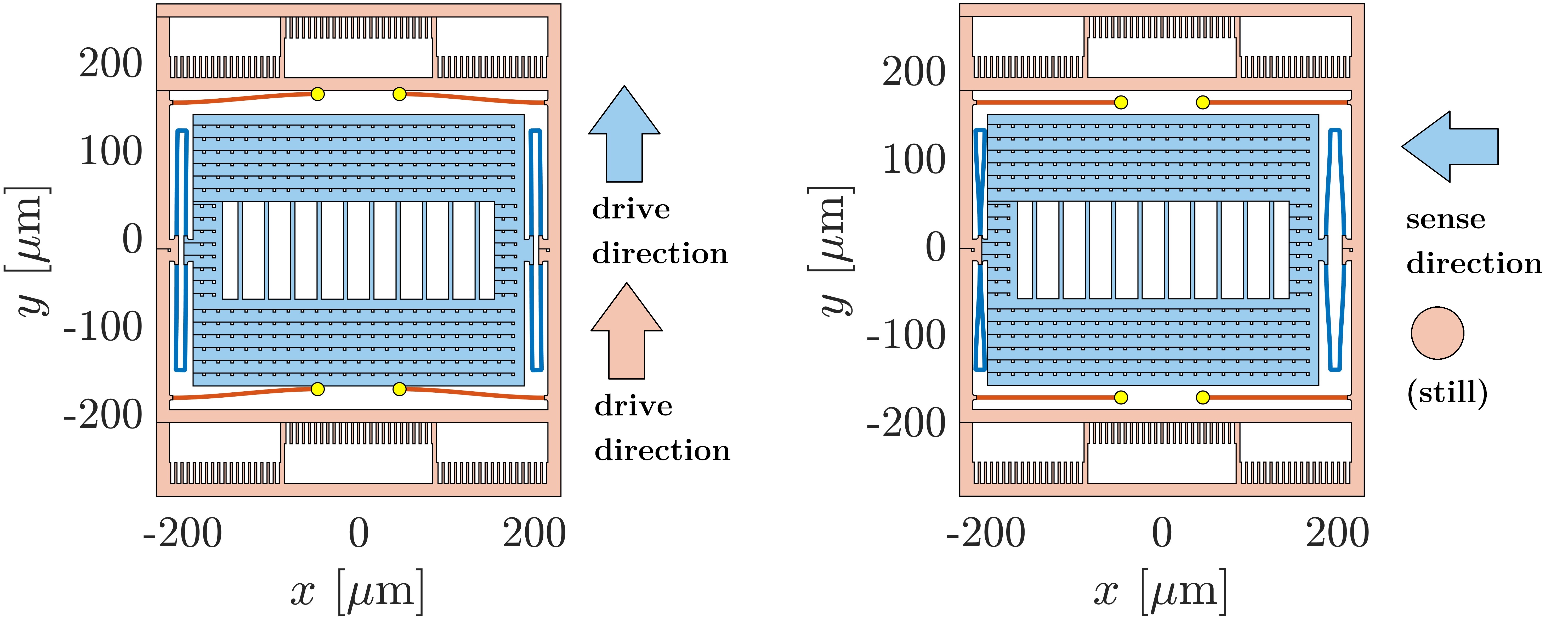}
    \caption{Drive (left) and sense (right) modes of the initial layout for the MEMS gyroscope. The folded beams (blue) are used to connect the drive frame (light red) to the sense mass (light blue), while the drive beams (red) connect the drive frame to the ground (yellow circles). Arrows show the displacements of the frame and the mass.}
    \label{fig:mems_initial}
\end{figure}

Consider the MEMS gyroscope prototype \cite{Marconi2021gyro} in Fig.~\ref{fig:mems_initial}. The drive (light red) and sense (light blue) frames are fixed throughout the optimization, as they host the actuation and sensing electrodes. Since the problem is 2D, each frame is modeled as a rigid mass with lumped interial properties and 3 degrees of freedom (DOFs) corresponding to the displacements and rotation of its center of mass. On the other hand, the drive (red) and sense (blue) beams are represented by geometrically nonlinear von K\'arm\'an beams. The material is polysilicon \cite{Acar2009}, characterized by Young's modulus $E = 148\,\mathrm{GPa}$, Poisson's ratio $\nu = 0.23$, and density $\rho = 2330\,\mathrm{kg/m^3}$. The Rayleigh damping parameters are $\alpha_R = 0.01$ and $\beta_R = 0$.

Using the Multi Point Constraint (MPC) approach \cite[Chapter~13]{cook2001concepts}, each mass (\textit{master}) is connected to a node (\textit{slave}) of the beam elements using the following \textit{linearized} kinematic constraint:
\begin{equation} \label{eq:mpc_constraint}
    \begin{bmatrix}
        u \\
        v \\
        \theta
    \end{bmatrix}_{slave}
    =
    \begin{bmatrix}
        1 \quad & 0 \quad & -\Delta y \\
        0 \quad & 1 \quad &  \Delta x \\
        0 \quad & 0 \quad & 1
    \end{bmatrix}
    \begin{bmatrix}
        u \\
        v \\
        \theta
    \end{bmatrix}_{master},
\end{equation}
being $u$, $v$, and $\theta$ the x-displacement, the y-displacement, and the rotation angle of a node, and where $\Delta x = x_{slave} - x_{master}$ (same for $\Delta y$).

The shape of the drive beams is parametrized using a superimposition of five harmonics according to
\begin{equation}
    y(x) = y_0(x) + \sum_{k = 1}^{5} A_k \sin\left(\frac{k\pi}{L}x\right),
\end{equation}
where $L$ is the length of the beam, $y_0(x)$ represents the original $y$ coordinates as if the beams were straight lines, and $A_k$ are the amplitudes of the harmonics. The other parameters of the optimization are the length and width of the drive beams ($L_d$ and $W_d$), the length and width of the sense beams ($L_s$ and $W_s$), and the width of the connecting elements between the sense beams $W_c$. As commonly done in MEMS devices \cite{Acar2009}, we impose a symmetry about the vertical axis, meaning that the two to drive beams, the two bottom drive beams, and the two folded beams are symmetric. Therefore, the geometry of the MEMS gyroscope is defined by a total number of 16 parameters.

\begin{table}[t!]
    \centering
    \begin{tabular}{l|lllll}
    \hline
    Sizes [$\mu$m]       &   $\Delta l_d = 4.21$ &   $\Delta l_s = 0.04$ &   $\Delta w_d = -0.39$ &   $\Delta w_s = -0.25$ &   $\Delta w_c = 0.04$ \\
    \hline
    Top amp. [$\mu$m]    & $A_1 = 0.00$ & $A_2 = 1.58$ & $A_3 =  0.02$ & $A_4 =  0.79$ & $A_5 = 0.06$ \\
    \hline
    Bottom amp. [$\mu$m] & $A_1 = 0.00$ & $A_2 = 1.58$ & $A_3 =  0.02$ & $A_4 =  0.79$ & $A_5 = 0.06$ \\
    \hline
    \end{tabular}
    \caption{Optimal design variables obtained solving problem~\eqref{opt:mems}. The top row lists the lengths and widths of the drive and sense beams, as variations with respect to the initial values. The middle and bottom rows show the amplitudes of the top and bottom drive beams, respectively.}
    \label{tab:mems_results}
\end{table}

The optimization problem aims at imposing the frequencies $\omega_{drive}$ and $\omega_{sense}$ corresponding to the drive and sense mode shapes, respectively (Fig.~\ref{fig:mems_initial}). In addition, a point $(\Omega_{drive}, x_{drive})$ is imposed for the backbone curve associated to the drive mode. To this end, the optimization problem is stated as
\begin{equation} \label{opt:mems}
    \begin{aligned}
        \min_{\bs{\mu}} \quad& J \\
        \mathrm{s.t.} \quad& \omega_{drive} - \omega_{drive,0} = 0 \\
        & \omega_{sense} - \omega_{sense,0} = 0 \\
        & \Omega_{drive} \left(\rho\left(x_{drive,0}\right)\right) - \Omega_{drive,0} = 0 \\
        & \bs{\mu}_L \leq \bs{\mu} \leq \bs{\mu}_U,
    \end{aligned}
\end{equation}
where the target values are $\omega_{drive,0} = 28\,\mathrm{kHz}$, $\omega_{sense,0} = 26\,\mathrm{kHz}$, $x_{drive,0} = 3\,\mathrm{\mu m}$, and $\Omega_{drive,0} = 26\,\mathrm{kHz}$. The beam lengths are allowed to vary between $-25\,\mathrm{\mu m}$ and  $25\,\mathrm{\mu m}$, while the beams widths between $-1\,\mathrm{\mu m}$ and $1\,\mathrm{\mu m}$. The drive beam amplitudes are bounded between $5\,\mathrm{\mu m}$ and $-5\,\mathrm{\mu m}$. At the first iteration, the design variables are all zero.

\begin{figure}[t!]\label{fig:mems_error}
    \centering
    \includegraphics[width=0.4\linewidth]{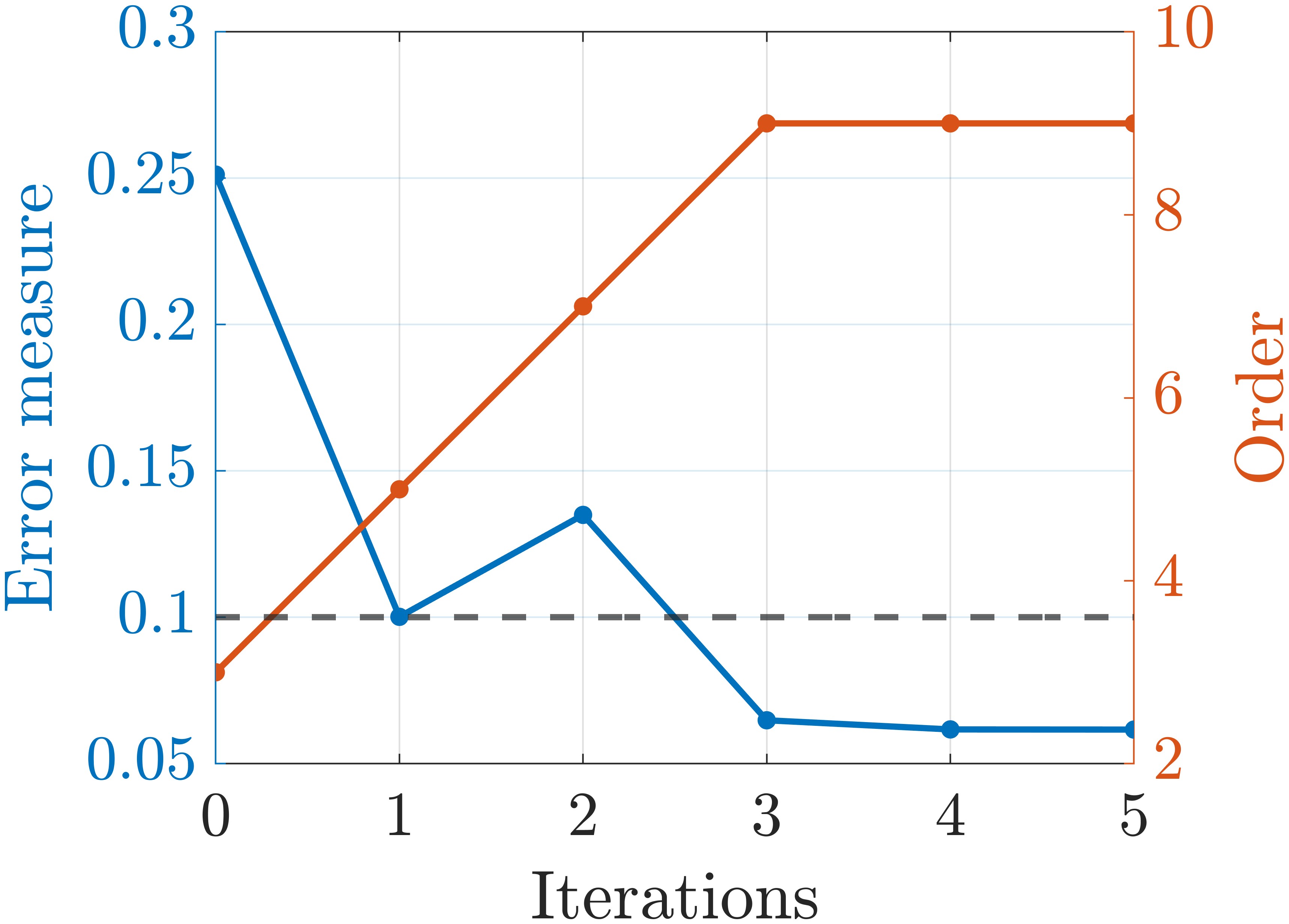}
    \caption{Evolution of the error measure $\varepsilon$ and expansion order during the optimization of the MEMS gyroscope.}
\end{figure}

Using the adjoint method for the sensitivity of the frequencies $\omega_{drive}$, $\omega_{sense}$, and $\Omega_{drive}$, the optimization reaches convergence in 5 iterations (around 2 minutes on a Windows laptop equipped with Intel Core i7-1255U CPU @ 1.70 GHz and 16.0 GB RAM). During the optimization, the expansion order is automatically adjusted to keep the error measure \cite{Pozzi2024backbone} below a threshold of $10^{-1}$. In particular, the expansion order has been increased from 3 to 9 (Fig.~\ref{fig:mems_error}).

\begin{figure}[t!]
    \centering
    \includegraphics[width=.9\linewidth]{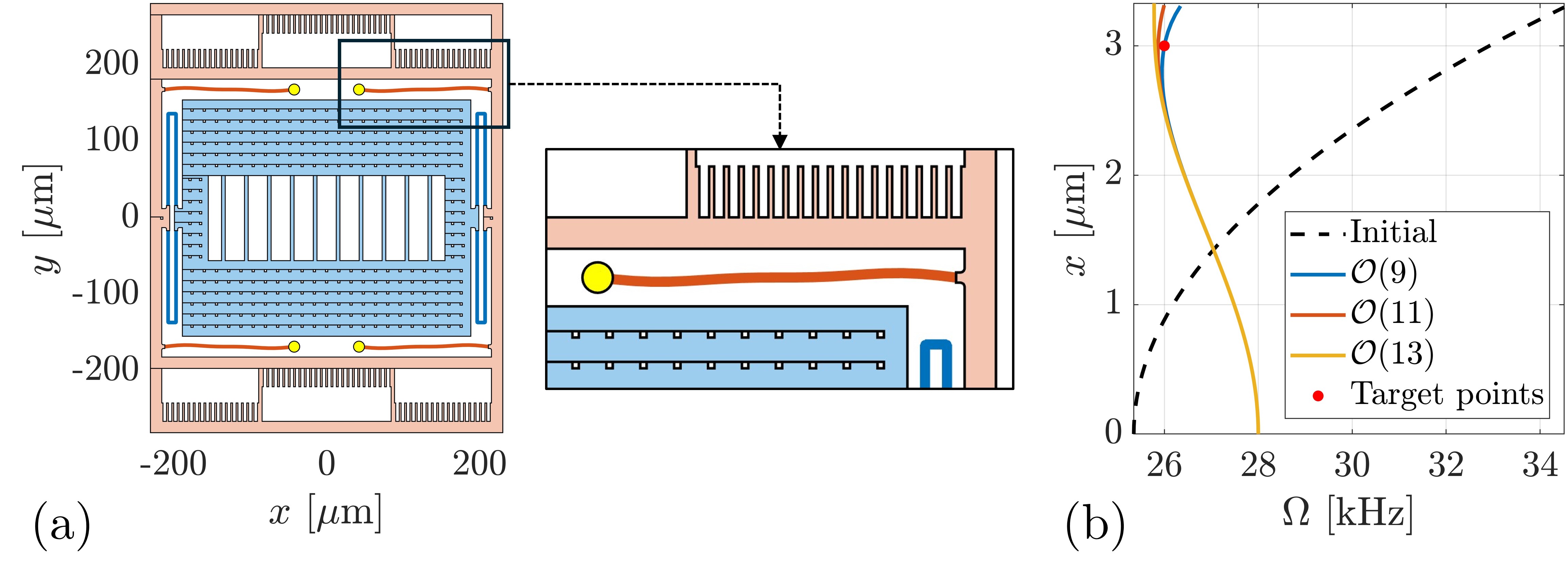}
    \caption{Solution of problem~\eqref{opt:mems}. (a) Optimal layout for the MEMS gyroscope. As it can be seen in the inset, the drive beams are not straight anymore. (b) Backbone curve associated to the drive mode of the optimized layout ($x$ is the RMS value of the drive frame y-displacement). Solutions at expansion orders 11 and 13 are reported for convergence analysis.}
    \label{fig:mems_solution_bb}
\end{figure}

The optimal layout is shown in Fig.~\ref{fig:mems_solution_bb}a, while the backbone curve associated with the drive mode is computed using two expansion orders (Fig.~\ref{fig:mems_solution_bb}b) to check the convergence of the SSM reduction.

\section{Conclusions} \label{sec:conclusion}

This work presents an optimization framework for tailoring the nonlinear dynamic response of lightly damped mechanical systems by leveraging Spectral Submanifold (SSM) reduction. The use of adjoint sensitivities in this formulation significantly reduces the computational cost compared to direct differentiation, making large-scale optimizations feasible. By formulating the backbone curve and its sensitivity directly in the physical space, the proposed method enables optimization of the physical amplitude-frequency relation.  The use of multi-index notation in our derivations overcomes the prior computational limitations of tensor-based formulation~\cite{Pozzi2024backbone}.  As our expressions enable optimization via SSMs up to arbitrary polynomial orders,  we demonstrated how the optimization accuracy can be systematically controlled by automatic adjustment of the SSM expansion order.

The results demonstrate the effectiveness of the proposed framework in optimizing the nonlinear response of mechanical systems. This approach broadens the scope of SSM-based methods, making them viable for practical engineering applications that require precise control over nonlinear dynamic behavior.

Future developments will focus on extending this approach to the optimization of forced-frequency responses, multiple modes with internal resonances, and parametric amplification.

%
%
%
%
%


\vspace{.5cm}
\noindent\textbf{Acknowledgments} M. Pozzi, J. Marconi and F. Braghin acknowledge the financial support of STMicroelectronics (award number 4000614871). Mingwu Li acknowledges the financial support of the National Natural Science Foundation of China (No. 12302014).

\newpage
\begin{appendices}

\section{On the Rayleigh damping coefficients} \label{app:damping}

The SSM formulation described hereby holds for lightly damped mechanical systems for which $\xi < 1$. By looking at Eq.~\eqref{eq:damping_ratio}, this condition becomes
\begin{equation} \label{eq:damping_condition}
    \alpha_R - 2 \omega + \beta_R \omega^2 < 0,
\end{equation}
where $\omega>0$. However, the eigenfrequency $\omega$ may change during the optimization process, and we cannot strictly ensure that $\xi < 1$ if $\alpha_R$ and $\beta_R$ are fixed. Assuming that $\alpha_R>0$ and $\beta_R>0$, and knowing that the roots of the equality associated to Eq.~\eqref{eq:damping_condition} are
\begin{equation}
    \omega_{1,2} = \frac{1 \pm \sqrt{1 - \alpha_R \beta_R}}{\beta_R},
\end{equation}
we have the following conditions on $\omega$ to satisfy Eq.~\eqref{eq:damping_condition}:
\begin{equation}
    \left\{
    \begin{aligned}
        \mathrm{if} &\quad \beta_R = 0 & \to &\quad \omega > \frac{\alpha_R}{2}\\
        \mathrm{if} &\quad \alpha_R = 0 & \to &\quad \omega < \frac{2}{\beta_R} \\
        \mathrm{if} &\quad 0 < \alpha_R\beta_R < 1 & \to &\quad \omega_1 < \omega < \omega_2 \\
        \mathrm{if} &\quad \alpha_R\beta_R > 1 & \to &\quad \text{Eq.~\eqref{eq:damping_condition} is never satisfied (i.e., $\xi> 1$)} \\
    \end{aligned}
    \right.
\end{equation}

In this work $\alpha_R$ and $\beta_R$ are selected and fixed before starting the optimization routine, so their values should be carefully chosen according to the conditions above, depending on the specific problem at hand. 

A possible alternative is to impose the value of $\xi$ and compute $\alpha_R$ and $\beta_R$ accordingly~\cite{Li2024adjoint}.

\section{Multi-index notation} \label{app:multi_index}

A multi-index $\bs{m} \in \mathbb{N}^M$ of order $m = |\bs{m}|_1$ is an $M$-dimensional vector for which addition, subtraction, and other operations are defined element-wise. For instance, the multi-indices can be used to define multivariate monomials of order $m$ as
\begin{equation}
    \bs{p}^{\bs{m}} = p_1^{m_1} \cdots p_M^{m_M}.
\end{equation}

To denote quantities related to a specific multi-index $\bs{m}=\{x,y\}$ we use the subscript $\bullet_{xy}$. For instance, the manifold coefficient related to $\bs{m}=\{2,1\}$ is $\bs{w}_{21}$.

Additionally, given a multi-index $\bs{m} = \{x, y\}$, define the symmetric multi-index $\bs{m}^s = \{y, x\}$.

\section{Near-resonance condition} \label{app:near_resonance}

Being $\bs{\Lambda}=\{\lambda_1,\lambda_2\}^T = \{\lambda,\bar{\lambda}\}^T$, we have
\begin{equation}
\begin{aligned}
    \Lambda_{\bs{m}} &= \bs{\Lambda}\cdot\bs{m} = m_1\lambda+m_2\bar\lambda =-(m_1+m_2)\alpha + \mathrm{i} (m_1 - m_2) \omega_d.
\end{aligned}
\end{equation}

For an \textit{undamped} system ($\xi = 0$), we have 
\begin{equation*}
    \Re(\lambda)=\Re(\bar\lambda)=0 \, \implies \, \Lambda_{\bs{m}} = \mathrm{i} (m_1-m_2)\omega.
\end{equation*}
In this case, an \textit{inner} resonance occurs if
\begin{equation}\label{eq:inner_res_condition}
\begin{aligned}
    \Lambda_{\bs{m}}&=\lambda_j \quad j=1,2\\
    \mathrm{i}(m_1-m_2)\omega&=\pm \mathrm{i}\omega \\
    m_1-m_2 &= \pm 1.
\end{aligned}
\end{equation}
As a consequence, we can conclude that at the \textit{even} order $m$ there are no inner resonances and the coefficients $\bs{R}_{\bs{m}}$ are null, whereas at the \textit{odd} order $m$ there are only 2 MIs such that
\begin{equation}
    \begin{cases}
        m_1 - m_2 = +1 \implies R^1_{\bs{m}} \neq 0 \\
        m_1 - m_2 = -1 \implies R^2_{\bs{m}} \neq 0.
    \end{cases}
\end{equation}
For a generic damped system instead, we write
\begin{equation}
    -(m_1+m_2)\alpha + \mathrm{i}(m_1-m_2)\omega_d = -\alpha \pm \mathrm{i}\omega_d,
\end{equation}
resulting into the following conditions:
\begin{equation}
    \begin{cases}
        m_1 - m_2 = \pm 1 \\
        m_1 + m_2 = 1.
    \end{cases}
\end{equation}
To strictly satisfy the above expressions, the only possible solutions are $\bs{m}=\{1,0\}$ and $\bs{m}=\{0,1\}$. However, under the assumption that $|\Re{\lambda}| \ll |\Im{\lambda}|$, we can safely ignore the second constraint in Eq.~\eqref{eq:inner_res_condition} and use the \textit{near-resonance} condition $\Lambda_{\bs{m}} \approx \lambda_j$.

\section{SSM computation} \label{app:ssm_steps}
\subsection{Leading order}

At the leading order, there are two multi-indices $\bs{m} = \{1, 0\}$ and $\bs{m} = \{0, 1\}$. The coefficients for the nonlinear dynamics are typically chosen as the eigenvalues of the system:
\begin{align} \label{eq:R1}
    R_{10}^{1} = \lambda_1 , \quad R_{01}^{2} = \lambda_2,
\end{align}
while $R_{10}^{2} = R_{01}^{1} = 0$. The manifold coefficients are then selected as
\begin{align} \label{eq:w1}
    \bs{w}_{10} = \bs{\phi} , \quad 
    \bs{w}_{01} = \bs{\phi}, \quad 
    \dot{\bs{w}}_{10} = \lambda_1 \bs{\phi} , \quad
    \dot{\bs{w}}_{01} = \lambda_2 \bs{\phi}.
\end{align}

%
%
%
%
%

\subsection{Higher orders}

After choosing the leading order coefficients, the higher order ones are computed iteratively starting with the lowest one. The first step is to compute the eigenvalue coefficient
\begin{equation} \label{eq:LambdaM}
    \Lambda_{\bs{m}} = \bs{m} \cdot \bs{\Lambda}.
\end{equation}

Then, using tensor notation, the nonlinear force contribution at the current order is
\begin{equation} \label{eq:fm}
    \begin{aligned}
        f^i_{\bs{m}} &= \sum_{\substack{\bs{u}, \bs{k} \in \mathbb{N}^2 \\ \bs{u} + \bs{k} = \bs{m}}} T_2^{ijk} w^j_{\bs{u}} w^k_{\bs{k}}  \quad + \sum_{\substack{\bs{u}, \bs{k}, \bs{l} \in \mathbb{N}^2 \\ \bs{u} + \bs{k} + \bs{l} = \bs{m}}} T_3^{ijkl} w^j_{\bs{u}} w^k_{\bs{k}} w^l_{\bs{l}}.
    \end{aligned}
\end{equation}

The vectors $\bs{V}_{\bs{m}}$ and $\dot{\bs{V}}_{\bs{m}}$ collect the contributions from the orders lower than the current one. These vectors are computed as
\begin{align} \label{eq:Vm}
    \bs{V}_{\bs{m}} = \sum_{j = 1}^2 \sum_{\substack{\bs{u}, \bs{k} \in \mathbb{N}^2 \\ 1 < u,k < m \\ \bs{u} + \bs{k} - \bs{e}_j = \bs{m}}} \bs{w}_{\bs{u}} u_j R_{\bs{k}}^j , \quad
    \dot{\bs{V}}_{\bs{m}} = \sum_{j = 1}^2 \sum_{\substack{\bs{u}, \bs{k} \in \mathbb{N}^2 \\ 1 < u,k < m \\ \bs{u} + \bs{k} - \bs{e}_j = \bs{m}}} \dot{\bs{w}}_{\bs{u}} u_j R_{\bs{k}}^j
\end{align}

The vectors $\bs{f}_{\bs{m}}$, $\bs{V}_{\bs{m}}$, and $\dot{\bs{V}}_{\bs{m}}$ are then combined into
\begin{equation} \label{eq:Cm}
    \bs{C}_{\bs{m}} = -\bb{M} \dot{\bs{V}}_{\bs{m}} - \left( \Lambda_{\bs{m}} \bb{M} + \bb{C} \right) \bs{V}_{\bs{m}} - \bs{f}_{\bs{m}}
\end{equation}

For multi-indices that are in near-resonance condition (Appendix~\ref{app:near_resonance}), the vector $\bs{C}_{\bs{m}}$ is used in the computation of the reduced dynamics coefficient:
\begin{equation} \label{eq:Rm}
    R_{\bs{m}}^j = \frac{\bs{\phi}^T \bs{C}_{\bs{m}}}{\Lambda_{\bs{m}} + \lambda_j + \alpha_R + \beta_R \omega^2}.
\end{equation}

Next, the left-hand side of the cohomological equation is assembled as
\begin{equation} \label{eq:Lm}
    \bb{L}_{\bs{m}} = \bb{K} + \Lambda_{\bs{m}} \bb{C} + \Lambda_{\bs{m}}^2 \bb{M},
\end{equation}
while the right-hand side reads
\begin{equation} \label{eq:hm}
    \bs{h}_{\bs{m}} = \bs{C}_{\bs{m}} + \sum_{j = 1}^2 \bs{D}_{\bs{m}}^j R_{\bs{m}}^j,
\end{equation}
where
\begin{equation} \label{eq:Dm}
    \bs{D}_{\bs{m}}^j = -\left[ \left( \Lambda_{\bs{m}} + \lambda_j \right) \bb{M} + \bb{C} \right] \bs{\phi}.
\end{equation}
Notice that $\bs{D}_{\bs{m}}^j$ is required only if $R_{\bs{m}}^j \neq 0$ (near-resonance condition).

The cohomological equation is solved for the manifold coefficient $\bs{w}_{\bs{m}}$:
\begin{equation} \label{eq:cohomological}
    \bb{L}_{\bs{m}} \bs{w}_{\bs{m}} = \bs{h}_{\bs{m}}
\end{equation}

Finally, the manifold velocity coefficient is computed as
\begin{equation} \label{eq:dwm}
    \dot{\bs{w}}_{\bs{m}} = \Lambda_{\bs{m}} \bs{w}_{\bs{m}} + \sum_{j = 1}^2 R_{\bs{m}}^j \bs{\phi} + \bs{V}_{\bs{m}}.
\end{equation}

\textbf{Property} \textit{(Symmetric multi-index)}. In a general setting, all these steps needs to be repeated for every multi-index at each order. However, for the type of problem under consideration, the coefficients associated with symmetric multi-indices are complex conjugate (e.g., $\bs{w}_{12} = \bar{\bs{w}}_{21}$). This property applies to all coefficients involved in the SSM computation and is leveraged to reduce the overall computational cost. Specifically, we only need to compute the coefficients corresponding to multi-indices $\bs{m} = \{x,y\}$ satisfying $x \geq y$. Then, the coefficients for the remaining multi-indices, where $x < y$, are obtained as the complex conjugates of their symmetric counterparts.\\

\section{Direct differentiation sensitivity} \label{app:direct_differentiation}

Using the chain rule, the derivative of the response frequency $\Omega$ (Eq.~\eqref{eq:Omega}) with respect to the design variable $\mu$ is
\begin{equation} \label{eq:Omega_dd}
    \begin{aligned}
        \dv{\Omega}{\mu} &= \frac{1}{2} \mathrm{i} \left( \dv{\lambda_2}{\mu} - \dv{\lambda_1}{\mu} \right) + \frac{1}{2} \mathrm{i} \sum_{m > 1} \left[ \left( \dv{R_{\bs{m}}^2}{\mu} - \dv{R_{\bs{m}}^1}{\mu} \right) \rho^{m - 1} +  \dv{\rho}{\mu} \left( R_{\bs{m}}^2 - R_{\bs{m}}^1 \right) (m - 1) \rho^{m - 2} \right].
    \end{aligned}
\end{equation}

The derivative of the reduced amplitude $\rho$ with respect to $\mu$ is computed by differentiating Eq.~\eqref{eq:x_rms}:
\begin{equation}
    \dv{\rho}{\mu} = -\frac{\sum_k x_k^i \sum_{\bs{m}} \dv{w_{\bs{m}}^i}{\mu} \bs{p}_k^{\bs{m}}}{\sum_k x_k^i \sum_{\bs{m}} w_{\bs{m}}^i m \rho^{m - 1} \tilde{\bs{p}}_k^{\bs{m}}},
\end{equation}
where $x_k^i = x^i(\theta_k)$, $\bs{p}_k = \bs{p}(\theta_k)$, and $\tilde{\bs{p}}_k = \{ e^{\mathrm{i} \theta_k}, e^{-\mathrm{i} \theta_k} \}^T$.

To evaluate the sensitivity of $\Omega$, the derivatives of $\lambda$, $R^j_{\bs{m}}$, and $\bs{w}_{\bs{m}}$ are required. Before computing them, we need to evaluate the derivatives of $\omega$ and $\bs{\phi}$ with respect to $\mu$. These are computed together by differentiating the generalized eigenvalue problem (Eq.~\eqref{eq:gep}) and the mass normalization condition (Eq.~\eqref{eq:mass_normalization}):
\begin{equation}
\begin{aligned}
    \begin{bmatrix}
        \bb{K} - \omega^2 \bb{M} & -2\omega \bb{M} \bs{\phi} \\
        -2\omega \bs{\phi}^T \bb{M} & 0
    \end{bmatrix}
    \begin{bmatrix}
        \dv{\bs{\phi}}{\mu} \\ \dv{\omega}{\mu}
    \end{bmatrix} = 
    \begin{bmatrix}
        \left( \omega^2 \dv{\bb{M}}{\mu} - \dv{\bb{K}}{\mu} \right) \bs{\phi} \\
        \omega \bs{\phi}^T \dv{\bb{M}}{\mu} \bs{\phi}.
    \end{bmatrix}
\end{aligned}
\end{equation}

Then, the derivative of the eigenvalue $\lambda$ is computed as
\begin{equation}
\begin{aligned}
    \dv{\lambda}{\mu} &= -\xi \dv{\omega}{\mu} - \omega \dv{\xi}{\mu} + \mathrm{i} \dv{\omega}{\mu} \sqrt{1 - \xi^2} - \mathrm{i} \omega \frac{\xi}{\sqrt{1 - \xi^2}} \dv{\xi}{\mu},
\end{aligned}
\end{equation}
where the derivative of the damping ratio $\xi$ is
\begin{equation} \label{eq:damping_ratio_dd}
    \dv{\xi}{\mu} = \frac{\beta_R \omega^2 - \alpha_R}{2 \omega^2} \dv{\omega}{\mu}.
\end{equation}

Next, the derivatives of the leading order SSM coefficients are
\begin{align}
    \dv{R_{10}^{1}}{\mu} &= \dv{\bar{R}_{01}^{2}}{\mu} = \dv{\lambda}{\mu}, \quad\quad \dv{\bs{w}_{10}}{\mu} = \dv{\bar{\bs{w}}_{01}}{\mu} = \dv{\bs{\phi}}{\mu}, \quad\quad \dv{\dot{\bs{w}}_{10}}{\mu} = \dv{\dot{\bar{\bs{w}}}_{01}}{\mu} = \dv{\lambda}{\mu} \bs{\phi} + \lambda\dv{\bs{\phi}}{\mu}.
\end{align}

After that, we can differentiate all the higher order SSM coefficients starting from the lowest one. In particular, the derivative of the eigenvalue coefficient is simply
\begin{equation}
    \dv{\Lambda_{\bs{m}}}{\mu} = \bs{m} \cdot \dv{\bs{\Lambda}}{\mu}.
\end{equation}

Using the tensor notation, it is possible to write the derivative of the nonlinear force contribution:
\begin{equation}
\begin{aligned}
    \dv{f_{\bs{m}}^i}{\mu} &= \sum_{\substack{\bs{u}, \bs{k} \in \mathbb{N}^2 \\ \bs{u} + \bs{k} = \bs{m}}} \left[ \dv{T_2^{ijk}}{\mu} w^j_{\bs{u}} w^k_{\bs{k}} + T_2^{ijk} \left( \dv{w^j_{\bs{u}}}{\mu} w^k_{\bs{k}} + w^j_{\bs{u}} \dv{w^k_{\bs{k}}}{\mu} \right) \right] \\
    &\quad + \sum_{\substack{\bs{u}, \bs{k}, \bs{l} \in \mathbb{N}^2 \\ \bs{u} + \bs{k} + \bs{l} = \bs{m}}} \left[ \dv{T_3^{ijkl}}{\mu} w^j_{\bs{u}} w^k_{\bs{k}} w^l_{\bs{l}} + T_3^{ijkl} \left( \dv{w^j_{\bs{u}}}{\mu} w^k_{\bs{k}} w^l_{\bs{l}} + w^j_{\bs{u}} \dv{w^k_{\bs{k}}}{\mu} w^l_{\bs{l}} + w^j_{\bs{u}} w^k_{\bs{k}} \dv{w^l_{\bs{l}}}{\mu} \right) \right].
\end{aligned}
\end{equation}

The derivatives of vectors $\bs{V}_{\bs{m}}$ and $\dot{\bs{V}}_{\bs{m}}$ are written as:
\begin{align}
    \dv{\bs{V}_{\bs{m}}}{\mu} &= \sum_{j = 1}^2 \sum_{\substack{\bs{u}, \bs{k} \in \mathbb{N}^2 \\ 1 < u,k < m \\ \bs{u} + \bs{k} - \bs{e}_j = \bs{m}}} u_j \left( \dv{\bs{w}_{\bs{u}}}{\mu} R_{\bs{k}}^j + \bs{w}_{\bs{u}} \dv{R_{\bs{k}}^j}{\mu} \right) \\
    \dv{\dot{\bs{V}}_{\bs{m}}}{\mu} &= \sum_{j = 1}^2 \sum_{\substack{\bs{u}, \bs{k} \in \mathbb{N}^2 \\ 1 < u,k < m \\ \bs{u} + \bs{k} - \bs{e}_j = \bs{m}}} u_j \left( \dv{\dot{\bs{w}}_{\bs{u}}}{\mu} R_{\bs{k}}^j + \dot{\bs{w}}_{\bs{u}} \dv{R_{\bs{k}}^j}{\mu} \right)
\end{align}

The derivative of vector $\bs{C}_{\bs{m}}$ reads
\begin{equation}
\begin{aligned}
    \dv{\bs{C}_{\bs{m}}}{\mu} &= -\dv{\bb{M}}{\mu} \dot{\bs{V}}_{\bs{m}} - \bb{M} \dv{\dot{\bs{V}}_{\bs{m}}}{\mu} - \left( \Lambda_{\bs{m}} \bb{M} + \bb{C} \right) \dv{\bs{V}_{\bs{m}}}{\mu} - \left( \dv{\Lambda_{\bs{m}}}{\mu} \bb{M} + \Lambda_{\bs{m}} \dv{\bb{M}}{\mu} + \dv{\bb{C}}{\mu} \right) \bs{V}_{\bs{m}} - \dv{\bs{f}_{\bs{m}}}{\mu}
\end{aligned}
\end{equation}

For the multi-indices that are in near-resonance condition (Appendix~\ref{app:near_resonance}), we also need to differentiate the reduced dynamics coefficient:
\begin{equation}
\begin{aligned}
    \dv{R_{\bs{m}}^j}{\mu} = \frac{\dv{\bs{\phi}^T}{\mu} \bs{C}_{\bs{m}} + \bs{\phi}^T \dv{\bs{C}_{\bs{m}}}{\mu}}{\Lambda_{\bs{m}} + \lambda_j + \alpha_R + \beta_R \omega^2} - R_{\bs{m}}^j \frac{\dv{\Lambda_{\bs{m}}}{\mu} + \dv{\lambda_j}{\mu} + 2\beta_R \omega \dv{\omega}{\mu}}{\Lambda_{\bs{m}} + \lambda_j + \alpha_R + \beta_R \omega^2}.
\end{aligned}
\end{equation}

Next, the derivative of the manifold coefficient is obtained by differentiating the cohomological equation:
\begin{equation}
    \bb{L}_{\bs{m}} \dv{\bs{w}_{\bs{m}}}{\mu} = \dv{\bs{h}_{\bs{m}}}{\mu} - \dv{\bb{L}_{\bs{m}}}{\mu} \bs{w}_{\bs{m}},
\end{equation}
where
\begin{align}
\begin{split}
    \dv{\bb{L}_{\bs{m}}}{\mu} &= \dv{\bb{K}}{\mu} + \Lambda_{\bs{m}} \dv{\bb{C}}{\mu} + \Lambda_{\bs{m}}^2 \dv{\bb{M}}{\mu} + \left( \bb{C} + 2\Lambda_{\bs{m}} \bb{M} \right) \dv{\Lambda_{\bs{m}}}{\mu}
\end{split} \\
    \dv{\bs{h}_{\bs{m}}}{\mu} &= \dv{\bs{C}_{\bs{m}}}{\mu} + \sum_{j = 1}^2 \left( \dv{\bs{D}_{\bs{m}}^j}{\mu} R_{\bs{m}}^j + \bs{D}_{\bs{m}}^j \dv{R_{\bs{m}}^j}{\mu} \right) \\
\begin{split}
    \dv{\bs{D}_{\bs{m}}^j}{\mu} &= -\left[ \left( \Lambda_{\bs{m}} + \lambda_j \right) \bb{M} + \bb{C} \right] \dv{\bs{\phi}}{\mu} - \left[ \left( \dv{\Lambda_{\bs{m}}}{\mu} + \dv{\lambda_j}{\mu} \right) \bb{M} + \left( \Lambda_{\bs{m}} + \lambda_j \right) \dv{\bb{M}}{\mu} + \dv{\bb{C}}{\mu} \right] \bs{\phi}.
\end{split}
\end{align}

Finally, the derivative of the manifold velocity coefficient is
\begin{equation}
\begin{aligned}
    \dv{\dot{\bs{w}}_{\bs{m}}}{\mu} &= \dv{\Lambda_{\bs{m}}}{\mu} \bs{w}_{\bs{m}} + \Lambda_{\bs{m}} \dv{\bs{w}_{\bs{m}}}{\mu} + \dv{\bs{V}_{\bs{m}}}{\mu} + \sum_{j = 1}^2 \left( \dv{R_{\bs{m}}^j}{\mu} \bs{\phi} + R_{\bs{m}}^j \dv{\bs{\phi}}{\mu} \right).
\end{aligned}
\end{equation}

\section{Adjoint sensitivity of the damped frequency} \label{app:omega_adjoint}

Using the adjoint method, and assuming mass normalized mode shapes (Eq.\eqref{eq:mass_normalization}), the sensitivity of the eigenfrequency $\omega$ with respect to the design variable $\mu$ is
\begin{equation}
    \dv{\omega}{\mu} = \frac{\bs{\phi}^T \left(\pdv{\bb{K}}{\mu} - \omega^2 \pdv{\bb{M}}{\mu} \right) \bs{\phi}}{2\omega}.
\end{equation}

Using this expression, it is possible to write the sensitivity of the damped frequency $\omega_d$ as
\begin{equation}
    \dv{\omega_d}{\mu} = \dv{\omega}{\mu} \sqrt{1 - \xi^2} - \omega \frac{\xi}{\sqrt{1 - \xi^2}} \dv{\xi}{\mu},
\end{equation}
where the derivative of the damping ratio $\xi$ is
\begin{equation}
    \dv{\xi}{\mu} = \frac{\beta_R \omega^2 - \alpha_R}{2 \omega^2} \dv{\omega}{\mu}.
\end{equation}

\section{Partial derivatives} \label{app:adjoint}
\subsection{Partial derivatives with respect to the reduced amplitude}

The partial derivative of the response frequency $\Omega$ with respect to the reduced amplitude $\rho$ is
\begin{equation} \label{eq:pOmega_rho}
    \pdv{\Omega}{\rho} = \frac{1}{2} \mathrm{i} \sum_{m > 1} \left( R_{\bs{m}}^2 - R_{\bs{m}}^1 \right) (m - 1) \rho^{m - 2}.
\end{equation}

The partial derivative of the RMS physical amplitude $x$ with respect to the reduced amplitude $\rho$ is
\begin{equation} \label{eq:pxi_rho}
    \pdv{x}{\rho} = \frac{1}{N_{\theta} x} \sum_k x_k^i \sum_{\bs{m}} w^i_{\bs{m}} m \rho^{m - 1} \tilde{\bs{p}}_k^{\bs{m}}
\end{equation}
where $x_k^i = x^i(\theta_k)$ and $\tilde{\bs{p}}_k = \{ e^{\mathrm{i} \theta_k}, e^{-\mathrm{i} \theta_k} \}^T$.\\

%
%
%
%
%

\subsection{Partial derivatives with respect to the manifold coefficients}

Before deep diving into the partial derivatives, it is important to highlight that quantities at a given order do not depend on the manifold coefficients at higher orders. Therefore, the partial derivative of a term at order $q$ with respect to the manifold coefficient at order $m$ is nonzero only if $q > m$.

Using tensor notation, the partial derivative of the nonlinear force contribution is
\begin{equation} \label{eq:pf_w}
    \begin{aligned}
        \pdv{f_{\bs{q}}^i}{w_{\bs{m}}^p} &= \sum_{\substack{\bs{u}, \bs{k} \in \mathbb{N}^2 \\ \bs{u} + \bs{k} = \bs{q}}} T_2^{ijk} \left( \pdv{w_{\bs{u}}^j}{w_{\bs{m}}^p} w_{\bs{k}}^k + w_{\bs{u}}^j \pdv{w_{\bs{k}}^k}{w_{\bs{m}}^p} \right) + \sum_{\substack{\bs{u}, \bs{k}, \bs{l} \in \mathbb{N}^2 \\ \bs{u} + \bs{k} + \bs{l} = \bs{q}}} T_3^{ijkl} \left( \pdv{w_{\bs{u}}^j}{w_{\bs{m}}^p} w_{\bs{k}}^k w_{\bs{l}}^l + w_{\bs{u}}^j \pdv{w_{\bs{k}}^k}{w_{\bs{m}}^p} w_{\bs{l}}^l + w_{\bs{u}}^j w_{\bs{k}}^k \pdv{w_{\bs{l}}^l}{w_{\bs{m}}^p} \right),
    \end{aligned}
\end{equation}
where
\begin{equation}
    \pdv{w_{\bs{u}}^j}{w_{\bs{m}}^p} = \left\{
    \begin{aligned}
        1 \quad &\mathrm{if} \quad \bs{u} = \bs{m} \quad \mathrm{and} \quad p = j \\
        0 \quad &\mathrm{otherwise}.
    \end{aligned}
    \right.
\end{equation}

The same condition applies for all the other partial derivatives in Eq.~\eqref{eq:pf_w}.

The partial derivative $\pdv{\bs{f}_{\bs{q}}}{\bs{w}_{\bs{m}}}$ is a matrix which, in general, can be dense. Therefore, to avoid storing the full matrix, the product $\bs{v}^T \pdv{\bs{f}_{\bs{q}}}{\bs{w}_{\bs{m}}}$ is stored instead, where $\bs{v}$ represents a vector that multiplies $\bs{f}_{\bs{q}}$ in the Lagrangian function.

The partial derivatives of vectors $\bs{V}_{\bs{q}}$ and $\dot{\bs{V}}_{\bs{q}}$ are
\begin{align} \label{eq:pV_w}
\begin{split}
    \bs{v}^T \pdv{\bs{V}_{\bs{q}}}{\bs{w}_{\bs{m}}} &= \sum_{j = 1}^2 \sum_{\substack{\bs{u}, \bs{k} \in \mathbb{N}^2 \\ 1 < u,k < q \\ \bs{u} + \bs{k} - \bs{e}_j = \bs{q}}} u_j \left( \bs{v}^T \pdv{\bs{w}_{\bs{u}}}{\bs{w}_{\bs{m}}} R_{\bs{k}}^j + (\bs{v}^T \bs{w}_{\bs{u}}) \pdv{R_{\bs{k}}^j}{\bs{w}_{\bs{m}}} \right)
\end{split} \\
\begin{split}
    \bs{v}^T \pdv{\dot{\bs{V}}_{\bs{q}}}{\bs{w}_{\bs{m}}} &= \sum_{j = 1}^2 \sum_{\substack{\bs{u}, \bs{k} \in \mathbb{N}^2 \\ 1 < u,k < q \\ \bs{u} + \bs{k} - \bs{e}_j = \bs{q}}} u_j \left( \bs{v}^T \pdv{\dot{\bs{w}}_{\bs{u}}}{\bs{w}_{\bs{m}}} R_{\bs{k}}^j + (\bs{v}^T \dot{\bs{w}}_{\bs{u}}) \pdv{R_{\bs{k}}^j}{\bs{w}_{\bs{m}}} \right),
\end{split}
\end{align}
where the products $\bs{v}^T \bs{w}_{\bs{u}}$ and $\bs{v}^T \dot{\bs{w}}_{\bs{u}}$ are scalar quantities. Moreover, the product $\bs{v}^T \pdv{\bs{w}_{\bs{u}}}{\bs{w}_{\bs{m}}}$ is equal to $\bs{v}^T$ if $\bs{u} = \bs{m}$, otherwise it is equal to a null row vector.

The partial derivatives of vector $\bs{C}_{\bs{q}}$ is
\begin{equation} \label{eq:pC_w}
\begin{aligned}
    \bs{v}^T \pdv{\bs{C}_{\bs{q}}}{\bs{w}_{\bs{m}}} &= -\bs{v}^T \bb{M} \pdv{\dot{\bs{V}}_{\bs{q}}}{\bs{w}_{\bs{m}}} - \bs{v}^T \pdv{\bs{f}_{\bs{q}}}{\bs{w}_{\bs{m}}} - \bs{v}^T \left( \Lambda_{\bs{q}} \bb{M} + \bb{C} \right) \pdv{\bs{V}_{\bs{q}}}{\bs{w}_{\bs{m}}}.
\end{aligned}
\end{equation}

The partial derivative of the reduced dynamics coefficient is
\begin{equation} \label{eq:pR_w}
    \pdv{R_{\bs{q}}^j}{\bs{w}_{\bs{m}}} = \frac{\bs{\phi}^T \pdv{\bs{C}_{\bs{q}}}{\bs{w}_{\bs{m}}}}{\Lambda_{\bs{q}} + \lambda_j + \alpha_R + \beta_R \omega^2},
\end{equation}
where $\bs{\phi}^T \pdv{\bs{C}_{\bs{q}}}{\bs{w}_{\bs{m}}}$ is computed using Eq.~\eqref{eq:pC_w} by setting $\bs{v} = \bs{\phi}$.

The partial derivative of the right-hand side of the cohomological equation is
\begin{equation} \label{eq:ph_w}
\begin{aligned}
    \bs{v}^T \pdv{\bs{h}_{\bs{q}}}{\bs{w}_{\bs{m}}} &= -\bs{v}^T \pdv{\bs{C}_{\bs{q}}}{\bs{w}_{\bs{m}}} - \sum_{j = 1}^2 (\bs{v}^T \bs{D}_{\bs{q}}^j) \pdv{R_{\bs{q}}^j}{\bs{w}_{\bs{m}}},
\end{aligned}
\end{equation}
where the product $\bs{v}^T \bs{D}_{\bs{q}}^j$ is a scalar quantity, and the vector $\bs{D}_{\bs{q}}^j$ does not depend on the manifold coefficients.

The partial derivative of the manifold velocity coefficient is
\begin{equation} \label{eq:pdw_w}
\begin{aligned}
    \bs{v}^T \pdv{\dot{\bs{w}}_{\bs{q}}}{\bs{w}_{\bs{m}}} &= \Lambda_{\bs{q}} \bs{v}^T \pdv{\bs{w}_{\bs{q}}}{\bs{w}_{\bs{m}}} + (\bs{v}^T \bs{\phi}) \sum_{j = 1}^2 \pdv{R_{\bs{q}}^j}{\bs{w}_{\bs{m}}} + \bs{v}^T \pdv{\bs{V}_{\bs{q}}}{\bs{w}_{\bs{m}}},
\end{aligned}
\end{equation}
where the product $\bs{v}^T \bs{\phi}$ is a scalar quantity.

The partial derivatives of the physical amplitude is
\begin{equation} \label{eq:pxi_wm}
    \pdv{x}{\bs{w}_{\bs{m}}} = \frac{1}{N_{\theta} x} \sum_k x_k^i \bs{a}^T \bs{p}_k^{\bs{m}},
\end{equation}
where $\bs{a}$ is a vector such that $a_j = 1$ if $i = j$, 0 otherwise.

The partial derivatives of the response frequency is
\begin{equation} \label{eq:pOmega_wm}
    \pdv{\Omega}{\bs{w}_{\bs{m}}} = \frac{1}{2} \mathrm{i} \sum_{q > m} \left( \pdv{R_{\bs{q}}^2}{\bs{w}_{\bs{m}}} - \pdv{R_{\bs{q}}^1}{\bs{w}_{\bs{m}}} \right) \rho^{m - 1}.
\end{equation}

%
%
%
%
%

\subsection{Partial derivatives with respect to the mode shape}

Using tensor notation, the partial derivative of the nonlinear force contribution is
\begin{equation} \label{eq:pf_phi}
    \begin{aligned}
        \pdv{f_{\bs{m}}^i}{\phi^p} &= \sum_{\substack{\bs{u}, \bs{k} \in \mathbb{N}^2 \\ \bs{u} + \bs{k} = \bs{m}}} T_2^{ijk} \left( \pdv{w_{\bs{u}}^j}{\phi^p} w_{\bs{k}}^k + w_{\bs{u}}^j \pdv{w_{\bs{k}}^k}{\phi^p} \right) + \sum_{\substack{\bs{u}, \bs{k}, \bs{l} \in \mathbb{N}^2 \\ \bs{u} + \bs{k} + \bs{l} = \bs{m}}} T_3^{ijkl} \left( \pdv{w_{\bs{u}}^j}{\phi^p} w_{\bs{k}}^k w_{\bs{l}}^l + w_{\bs{u}}^j \pdv{w_{\bs{k}}^k}{\phi^p} w_{\bs{l}}^l + w_{\bs{u}}^j w_{\bs{k}}^k \pdv{w_{\bs{l}}^l}{\phi^p} \right),
    \end{aligned}
\end{equation}
where the partial derivative $\pdv{w_{\bs{u}}^j}{\phi^p}$ is defined as follow:
\begin{equation}
    \pdv{w_{\bs{u}}^j}{\phi^p} = \left\{
    \begin{aligned}
        1 \quad &\mathrm{if} \quad |\bs{u}|_1 = 1 \quad \mathrm{and} \quad p = j \\
        0 \quad &\mathrm{otherwise}
    \end{aligned}
    \right.
\end{equation}

The same applies for all the other partial derivatives in Eq.~\eqref{eq:pf_phi}. As before, the product $\bs{v}^T \pdv{\bs{f}_{\bs{m}}}{\bs{\phi}}$ is computed and stored instead of the full derivative matrix.

The partial derivatives of vectors $\bs{V}_{\bs{m}}$ and $\dot{\bs{V}}_{\bs{m}}$ are
\begin{align} \label{eq:pV_phi}
\begin{split}
    \bs{v}^T \pdv{\bs{V}_{\bs{m}}}{\bs{\phi}} &= \sum_{j = 1}^2 \sum_{\substack{\bs{u}, \bs{k} \in \mathbb{N}^2 \\ 1 < u,k < m \\ \bs{u} + \bs{k} - \bs{e}_j = \bs{m}}} u_j (\bs{v}^T \bs{w}_{\bs{u}}) \pdv{R_{\bs{k}}^j}{\bs{\phi}}
\end{split} \\
\begin{split}
    \bs{v}^T \pdv{\dot{\bs{V}}_{\bs{m}}}{\bs{\phi}} &= \sum_{j = 1}^2 \sum_{\substack{\bs{u}, \bs{k} \in \mathbb{N}^2 \\ 1 < u,k < m \\ \bs{u} + \bs{k} - \bs{e}_j = \bs{m}}} u_j \left( \bs{v}^T \pdv{\dot{\bs{w}}_{\bs{u}}}{\bs{\phi}} R_{\bs{k}}^j + (\bs{v}^T \dot{\bs{w}}_{\bs{u}}) \pdv{R_{\bs{k}}^j}{\bs{\phi}} \right),
\end{split}
\end{align}
where the products $\bs{v}^T \bs{w}_{\bs{u}}$ and $\bs{v}^T \dot{\bs{w}}_{\bs{u}}$ are scalar quantities.

The partial derivatives of vector $\bs{C}_{\bs{m}}$ is
\begin{equation} \label{eq:pC_phi}
\begin{aligned}
    \bs{v}^T \pdv{\bs{C}_{\bs{m}}}{\bs{\phi}} &= -\bs{v}^T \bb{M} \pdv{\dot{\bs{V}}_{\bs{m}}}{\bs{\phi}} - \bs{v}^T \pdv{\bs{f}_{\bs{m}}}{\bs{\phi}} - \bs{v}^T \left( \Lambda_{\bs{m}} \bb{M} + \bb{C} \right) \pdv{\bs{V}_{\bs{q}}}{\bs{\phi}}.
\end{aligned}
\end{equation}

The partial derivative of the reduced dynamics coefficient is
\begin{equation} \label{eq:pR_phi}
    \pdv{R_{\bs{m}}^j}{\bs{\phi}} = \frac{\bs{C}_{\bs{m}}^T + \bs{\phi}^T \pdv{\bs{C}_{\bs{m}}}{\bs{\phi}}}{\Lambda_{\bs{q}} + \lambda_j + \alpha_R + \beta_R \omega^2}
\end{equation}
where $\bs{\phi}^T \pdv{\bs{C}_{\bs{m}}}{\bs{\phi}}$ is computed using Eq.~\eqref{eq:pC_phi} by setting $\bs{v} = \bs{\phi}$.

The partial derivatives of vector $\bs{D}_{\bs{m}}^j$ is
\begin{equation} \label{eq:pD_phi}
    \bs{v}^T \pdv{\bs{D}_{\bs{m}}^j}{\bs{\phi}} = - \bs{v}^T \left[ \left( \Lambda_{\bs{m}} + \lambda_j \right) \bb{M} + \bb{C} \right].
\end{equation}

The partial derivative of the right-hand side of the cohomological equation is
\begin{equation} \label{eq:ph_phi}
\begin{aligned}
    \bs{v}^T \pdv{\bs{h}_{\bs{m}}}{\bs{\phi}} &= -\bs{v}^T \pdv{\bs{C}_{\bs{m}}}{\bs{\phi}} - \sum_{j = 1}^2 (\bs{v}^T \bs{D}_{\bs{m}}^j) \pdv{R_{\bs{m}}^j}{\bs{\phi}} - \sum_{j = 1}^2 \bs{v}^T \pdv{\bs{D}_{\bs{m}}^j}{\bs{\phi}} R_{\bs{m}}^j
\end{aligned}
\end{equation}
where the product $\bs{v}^T \bs{D}_{\bs{m}}^j$ is a scalar quantity.

The partial derivative of the manifold velocity coefficient is
\begin{equation} \label{eq:pdw_phi}
\begin{aligned}
    \bs{v}^T \pdv{\dot{\bs{w}}_{\bs{m}}}{\bs{\phi}} &= (\bs{v}^T \bs{\phi}) \sum_{j = 1}^2 \pdv{R_{\bs{m}}^j}{\bs{\phi}} + \bs{v}^T \sum_{j = 1}^2 R_{\bs{m}}^j + \bs{v}^T \pdv{\bs{V}_{\bs{m}}}{\bs{\phi}}
\end{aligned}
\end{equation}
where the product $\bs{v}^T \bs{\phi}$ is a scalar quantity.

The partial derivatives of the physical amplitude is
\begin{equation} \label{eq:pxi_phi}
    \pdv{x}{\bs{\phi}} = \frac{1}{N_{\theta} x} \sum_k x_k^i \bs{a}^T (\bs{p}_k^{10} + \bs{p}_k^{01})
\end{equation}
where $\bs{a}$ is a vector such that $a_j = 1$ if $i = j$, 0 otherwise.

The partial derivatives of the response frequency is
\begin{equation} \label{eq:pOmega_phi}
    \pdv{\Omega}{\bs{\phi}} = \frac{1}{2} \mathrm{i} \sum_{m > 1} \left( \pdv{R_{\bs{m}}^2}{\bs{\phi}} - \pdv{R_{\bs{m}}^1}{\bs{\phi}} \right) \rho^{m - 1}
\end{equation}

%
%
%
%
%

\subsection{Partial derivatives with respect to the natural frequency}

Here, the partial derivatives are taken with respect to the natural frequency $\omega$, which is a scalar. Therefore, there is no need to multiply the partial derivatives by $\bs{v}^T$ as done in the previous cases.

The partial derivative of the eigenvalue coefficient is
\begin{equation} \label{eq:pLambdaM_omega}
    \pdv{\Lambda_{\bs{m}}}{\omega} = \bs{m} \cdot \pdv{\Lambda}{\omega},
\end{equation}
where
\begin{align}
    \pdv{\Lambda}{\omega} &=
    \begin{bmatrix}
        \pdv{\lambda}{\omega} \\ \pdv{\bar{\lambda}}{\omega}
    \end{bmatrix} \\
    \begin{split}
        \pdv{\lambda}{\omega} &= -\xi - \omega \pdv{\xi}{\omega} + \mathrm{i} \sqrt{1 - \xi^2} - \mathrm{i} \omega \frac{\xi}{\sqrt{1 - \xi^2}} \pdv{\xi}{\omega}
    \end{split} \\
    \dv{\xi}{\mu} &= \frac{\beta_R \omega^2 - \alpha_R}{2 \omega^2}.
\end{align}

The nonlinear force contribution $\bs{f}_{\bs{m}}$ does not depend on the eigenfrequency $\omega$, thus $\pdv{\bs{f}_{\bs{m}}}{\omega} = \bs{0}$.

The partial derivatives of vectors $\bs{V}_{\bs{m}}$ and $\dot{\bs{V}}_{\bs{m}}$ are
\begin{align} \label{eq:pV_omega}
    \pdv{\bs{V}_{\bs{m}}}{\omega} &= \sum_{j = 1}^2 \sum_{\substack{\bs{u}, \bs{k} \in \mathbb{N}^2 \\ 1 < u,k < m \\ \bs{u} + \bs{k} - \bs{e}_j = \bs{m}}} u_j \bs{w}_{\bs{u}} \pdv{R_{\bs{k}}^j}{\omega} \\
    \begin{split}
        \pdv{\dot{\bs{V}}_{\bs{m}}}{\omega} &= \sum_{j = 1}^2 \sum_{\substack{\bs{u}, \bs{k} \in \mathbb{N}^2 \\ 1 < u,k < m \\ \bs{u} + \bs{k} - \bs{e}_j = \bs{m}}} u_j \left( \pdv{\dot{\bs{w}}_{\bs{u}}}{\omega} R_{\bs{k}}^j + \dot{\bs{w}}_{\bs{u}} \pdv{R_{\bs{k}}^j}{\omega} \right).
    \end{split}
\end{align}

The partial derivatives of vector $\bs{C}_{\bs{m}}$ is
\begin{equation} \label{eq:pC_omega}
\begin{aligned}
    \pdv{\bs{C}_{\bs{m}}}{\omega} &= - \bb{M} \pdv{\dot{\bs{V}}_{\bs{m}}}{\omega} - \pdv{\bs{\Lambda}_{\bs{m}}}{\omega} \bb{M} \bs{V}_{\bs{m}}  - \left( \Lambda_{\bs{m}} \bb{M} + \bb{C} \right) \pdv{\bs{V}_{\bs{m}}}{\omega}.
\end{aligned}
\end{equation}

The partial derivative of the reduced dynamics coefficient is
\begin{equation} \label{eq:pR_omega}
\begin{aligned}
    \pdv{R_{\bs{m}}^j}{\omega} &= \frac{\bs{\phi}^T \pdv{\bs{C}_{\bs{m}}}{\omega}}{\Lambda_{\bs{m}} + \lambda_j + \alpha_R + \beta_R \omega^2} - R_{\bs{m}}^j \frac{\pdv{\Lambda_{\bs{m}}}{\omega} + \pdv{\lambda_j}{\omega} + 2\beta_R \omega}{\Lambda_{\bs{m}} + \lambda_j + \alpha_R + \beta_R \omega^2}.
\end{aligned}
\end{equation}

The partial derivative of the right-hand side of the cohomological equation is
\begin{equation} \label{eq:ph_omega}
    \pdv{\bs{h}_{\bs{m}}}{\omega} = \pdv{\bs{C}_{\bs{m}}}{\omega} + \sum_{j = 1}^2 \left( \pdv{\bs{D}_{\bs{m}}^j}{\omega} R_{\bs{m}}^j + \bs{D}_{\bs{m}}^j \pdv{R_{\bs{m}}^j}{\omega} \right).
\end{equation}
where
\begin{equation} \label{eq:pD_omega}
    \pdv{\bs{D}_{\bs{m}}^j}{\omega} = - \left( \pdv{\Lambda_{\bs{m}}}{\omega} + \pdv{\lambda_j}{\omega} \right) \bb{M} \bs{\phi}.
\end{equation}

The partial derivative of the left-hand side of the cohomological equation is
\begin{equation} \label{eq:pL_omega}
    \pdv{\bb{L}_{\bs{m}}}{\omega} = \left( \bb{C} + 2\Lambda_{\bs{m}} \bb{M} \right) \pdv{\Lambda_{\bs{m}}}{\omega}.
\end{equation}

The partial derivative of the manifold velocity coefficient is
\begin{equation} \label{eq:pdw_omega}
    \pdv{\dot{\bs{w}}_{\bs{m}}}{\omega} = \pdv{\Lambda_{\bs{m}}}{\omega} \bs{w}_{\bs{m}} + \bs{\phi} \sum_{j = 1}^2 \pdv{R_{\bs{m}}^j}{\omega} + \pdv{\bs{V}_{\bs{m}}}{\omega}.
\end{equation}

The RMS physical amplitude $x$ does not directly depend on $\omega$, while the partial derivative of the response frequency $\Omega$ is
\begin{equation} \label{eq:pOmega_omega}
\begin{aligned}
    \pdv{\Omega}{\omega} &= \frac{1}{2} \mathrm{i} \left( \pdv{\bar{\lambda}}{\omega} - \pdv{\lambda}{\omega} \right) + \frac{1}{2} \mathrm{i} \sum_{m > 1} \left( \pdv{R_{\bs{m}}^2}{\omega} - \pdv{R_{\bs{m}}^1}{\omega} \right) \rho^{m - 1}.
\end{aligned}
\end{equation}

%
%
%
%
%

\subsection{Partial derivatives with respect to the design variables}

Using tensor notation, the partial derivative of the nonlinear force contribution is
\begin{equation} \label{eq:pf_mu}
    \begin{aligned}
        \pdv{f^i_{\bs{m}}}{\mu} &= \sum_{\substack{\bs{u}, \bs{k} \in \mathbb{N}^2 \\ \bs{u} + \bs{k} = \bs{m}}} \pdv{T_2^{ijk}}{\mu} w^j_{\bs{u}} w^k_{\bs{k}} \quad + \sum_{\substack{\bs{u}, \bs{k}, \bs{l} \in \mathbb{N}^2 \\ \bs{u} + \bs{k} + \bs{l} = \bs{m}}} \pdv{T_3^{ijkl}}{\mu} w^j_{\bs{u}} w^k_{\bs{k}} w^l_{\bs{l}}.
    \end{aligned}
\end{equation}

The partial derivatives of vectors $\bs{V}_{\bs{m}}$ and $\dot{\bs{V}}_{\bs{m}}$ are
\begin{align} \label{eq:pV_mu}
    \pdv{\bs{V}_{\bs{m}}}{\mu} &= \sum_{j = 1}^2 \sum_{\substack{\bs{u}, \bs{k} \in \mathbb{N}^2 \\ 1 < u,k < m \\ \bs{u} + \bs{k} - \bs{e}_j = \bs{m}}} \bs{w}_{\bs{u}} u_j \pdv{R_{\bs{k}}^j}{\mu} \\
    \begin{split}
        \pdv{\dot{\bs{V}}_{\bs{m}}}{\mu} &= \sum_{j = 1}^2 \sum_{\substack{\bs{u}, \bs{k} \in \mathbb{N}^2 \\ 1 < u,k < m \\ \bs{u} + \bs{k} - \bs{e}_j = \bs{m}}} u_j \left[ \pdv{\dot{\bs{w}}_{\bs{u}}}{\mu} R_{\bs{k}}^j + \dot{\bs{w}}_{\bs{u}} \pdv{R_{\bs{k}}^j}{\mu} \right].
    \end{split}
\end{align}

The partial derivatives of vector $\bs{C}_{\bs{m}}$ is
\begin{equation} \label{eq:pC_mu}
\begin{aligned}
    \pdv{\bs{C}_{\bs{m}}}{\mu} &= -\pdv{\bb{M}}{\mu} \dot{\bs{V}}_{\bs{m}} - \bb{M} \pdv{\dot{\bs{V}}_{\bs{m}}}{\mu} - \left( \Lambda_{\bs{m}} \bb{M} + \bb{C} \right) \pdv{\bs{V}_{\bs{m}}}{\mu} - \left( \Lambda_{\bs{m}} \pdv{\bb{M}}{\mu} + \pdv{\bb{C}}{\mu} \right) \bs{V}_{\bs{m}} - \pdv{\bs{f}_{\bs{m}}}{\mu}.
\end{aligned}
\end{equation}

The partial derivative of the reduced dynamics coefficient is
\begin{equation} \label{eq:pR_mu}
    \pdv{R_{\bs{m}}^j}{\mu} = \frac{\bs{\phi}^T \pdv{\bs{C}_{\bs{m}}}{\mu}}{\Lambda_{\bs{m}} + \lambda_j + \alpha_R + \beta_R \omega^2}.
\end{equation}

The partial derivative of the right-hand side of the cohomological equation is
\begin{equation} \label{eq:ph_mu}
    \pdv{\bs{h}_{\bs{m}}}{\mu} = \pdv{\bs{C}_{\bs{m}}}{\mu} + \sum_{j = 1}^2 \left( \pdv{\bs{D}_{\bs{m}}^j}{\mu} R_{\bs{m}}^j + \bs{D}_{\bs{m}}^j \pdv{R_{\bs{m}}^j}{\mu} \right),
\end{equation}
where
\begin{equation} \label{eq:pD_mu}
    \pdv{\bs{D}_{\bs{m}}^j}{\mu} = -\left[ \left( \Lambda_{\bs{m}} + \lambda_j \right) \pdv{\bb{M}}{\mu} + \pdv{\bb{C}}{\mu} \right] \bs{\phi}.
\end{equation}

The partial derivative of the left-hand side of the cohomological equation is
\begin{equation} \label{eq:pL_mu}
    \pdv{\bb{L}_{\bs{m}}}{\mu} = \pdv{\bb{K}}{\mu} + \Lambda_{\bs{m}} \pdv{\bb{C}}{\mu} + \Lambda_{\bs{m}}^2 \pdv{\bb{M}}{\mu}.
\end{equation}

The partial derivative of the manifold velocity coefficient is
\begin{equation} \label{eq:pdw_mu}
    \pdv{\dot{\bs{w}}_{\bs{m}}}{\mu} = \sum_{j = 1}^2 \pdv{R_{\bs{m}}^j}{\mu} \bs{\phi} + \pdv{\bs{V}_{\bs{m}}}{\mu}.
\end{equation}

The RMS physical amplitude $x$ does not directly depend on $\mu$, while the partial derivative of the response frequency $\Omega$ is
\begin{equation} \label{eq:pOmega_mu}
    \pdv{\Omega}{\mu} = \frac{1}{2} \mathrm{i} \sum_{m > 1} \left( \pdv{R_{\bs{m}}^2}{\mu} - \pdv{R_{\bs{m}}^1}{\mu} \right) \rho^{m - 1}.
\end{equation}

\end{appendices}

\bibliographystyle{unsrt}
\bibliography{bibliography}

\end{document}